\documentclass[a4paper, secthm, amsthm]{elsart}
\usepackage{amsmath,amsfonts,euscript} 
\def\dis{\displaystyle}         
\def\eps{\varepsilon}
\def\Tb{\overline{T}}
\def\vv{\underline{v}}
\def\ww{\underline{w}}
\def\ff{\underline{f}}
\def\Tt{\tilde{T}}
\def\Lt{\tilde{L}}
\def\Ct{\tilde{C}}
\def\vt{\tilde{v}}
\def\ut{\tilde{u}}
\def\gt{\tilde{g}}
\def\at{\tilde{\alpha}}
\def\At{\tilde{A}}
\def\ldk{\lambda_{k}}
\def\lkp{\lambda_{k+1}}

\def\ldn{\lambda_{n}}
\def\lnp{\lambda_{n+1}}
\def\lnm{\lambda_{n-1}}
\def\sldn{\sqrt{\lambda_{n}}}
\def\slnp{\sqrt{\lambda_{n+1}}}
\def\slnm{\sqrt{\lambda_{n-1}}}
\def\sldo{\sqrt{\lambda_{1}}}

\def\R{{\mathbb R}}

\def\adh#1{\, \overline{#1}}

\def\sinc{\mathop{\rm sinc}}

\def\dist{\mathop{\rm dist}}
\newcommand{\cd}[1]{\mathop{\bf #1}}

\newcommand{\res}[1]{_{\rceil {\normalsize #1}}}

\begin{document}
%\pagestyle{headings}
%\markboth{Null-controllability cost for the heat equation.}{\sc L.Miller} 
%\markright{Null-controllability cost for the heat equation.}

%%%%%%%%%%%%%%%%%%%%%%%%%%%%%%%%%%%%%%%%%%%%%%%%%%%%%%%%%%%%
%% FRONT
%%%%%%%%%%%%%%%%%%%%%%%%%%%%%%%%%%%%%%%%%%%%%%%%%%%%%%%%%%%%

\begin{frontmatter}
\title{
%[null-controllability cost for the heat equation]
Geometric bounds on \\ 
the growth rate of null-controllability cost \\
for the heat equation in small time}
\author{Luc Miller}
\address{
{\'E}quipe Modal'X, JE 421, \\
Universit{\'e} Paris X, B{\^a}t.~G,
200 Av.~de la R{\'e}publique,
92001 Nanterre, France.}
\address{Centre de Math{\'e}matiques, UMR CNRS 7640, \\ 
{\'E}cole Polytechnique, 91128
Palaiseau, France.}
\ead{miller@math.polytechnique.fr}

\thanks{This work was partially supported by the ACI grant
  ``{\'E}quation des ondes : oscillations, dispersion et contr{\^o}le''.}
\date{July 11, 2003}%\date{\today}

\begin{abstract}
Given a control region $\Omega$ on a compact Riemannian manifold $M$,
we consider the heat equation with a source term $g$ localized in $\Omega$.
It is known that any initial data in $L^{2}(M)$ can be steered to $0$ 
in an arbitrarily small time $T$ 
by applying a suitable control $g$ in $L^{2}([0,T]\times\Omega)$,
and, as $T$ tends to $0$, 
the norm of $g$ grows like $\exp(C/T)$ times the norm of the data.
We investigate how $C$ depends on the geometry of $\Omega$. %% 72 words
We prove $C\geq d^{2}/4$ 
where $d$ is the largest distance of a point in $M$ from $\Omega$.
When $M$ is a segment of length $L$ controlled at one end, 
we prove $C\leq \alpha_{*}L^{2}$ for some $\alpha_{*}<2$.
Moreover, this bound implies $C\leq\alpha_{*}L_{\Omega}^{2}$
where $L_{\Omega}$ is the length of the longest generalized geodesic in $M$
which does not intersect $\Omega$.
The {\em control transmutation method} 
used in proving this last result 
is of a broader interest.
\end{abstract}


%% \begin{keyword}
%% Heat equation \sep control cost \sep
%% null-controllability \sep observabillity \sep  
%% small time asymptotics \sep 
%% multipliers \sep entire functions \sep transmutation.
%% \end{keyword}

%\thanks{2000 {\it Mathematics Subject Classification. } }
%35B37 PDE in connection with control problems, 35K05 Heat equation, 
%58J35 Heat and other parabolic equation methods
%58Jxx Partial differential equations on manifolds; differential operators

\end{frontmatter}

%%%%%%%%%%%%%%%%%%%%%%%%%%%%%%%%%%%%%%%%%%%%%%%%%%%%%%%%%%%%
%% THE PROBLEM
%%%%%%%%%%%%%%%%%%%%%%%%%%%%%%%%%%%%%%%%%%%%%%%%%%%%%%%%%%%%

\section{The problem}
\label{sec:pb}

Let $(M,g)$ be a smooth connected compact %oriented
$n$-dimensional Riemannian manifold with metric $g$ 
and boundary $\partial M$. 
When $\partial M\neq\emptyset$, $M$ denotes the interior 
and $\adh{M}=M\cup\partial{M}$.
Let $\dist:\adh{M}^{2}\to \mathbb{R}_{+}$ 
denote the distance function.
Let $\Delta$ denote the (negative) Dirichlet Laplacian on $L^2(M)$
with domain $D(\Delta)=H^1_0(M) \cap H^2(M)$.
%Let $\nu$ denote the exterior normal vector field.
%The Sobolev spaces and the energy (\ref{eqE}) are defined 
%with respect to the measure $dx_{g}$. 
%In local coordinates $(x^{1},x^{2},\dots,x^{n})$~:
%$dx_{g}=\sqrt{\det g}dx^{1}\cdots dx^{n}$.

Consider a positive control time $T$, 
and an open control region $\Omega$. %$\subset]0,T[\times M$ 
%such that \linebreak $\adh{\Omega}\subset M$. 
Let ${\bf 1}_{]0,T[\times \Omega}$ 
denote the characteristic function of the space-time 
control region $]0,T[\times \Omega$.
The heat equation on $M$ is said to be 
{\em null-controllable}
(or exactly controllable to zero) 
in time $T$ by interior controls on $\Omega$ 
if for all  
$u_{0}\in L^{2}(M)$
there is a control function 
$g\in L^{2}(\mathbb{R}\times M)$
such that the solution $u\in C^{0}([0,\infty),L^{2}(M))$
of the mixed Dirichlet-Cauchy problem:
\begin{equation} 
\label{eqHeat}
\partial_{t}u - \Delta u=\cd{1}_{]0,T[\times \Omega} g  
\quad {\rm in}\ ]0,T[\times M, \quad 
u=0 \quad {\rm on}\ ]0,T[\times\partial M,
\end{equation}
with Cauchy data
$u=u_{0}$ at $t=0$,
satisfies 
$u=0$ at $t=T$.
For a survey on this problem prior to 1978 we refer to \cite{Rus78}.
For a recent update, we refer to \cite{Zua01}.
Lebeau and Robbiano have proved 
(in \cite{LR95} using local Carleman estimates) 
that there is a continuous linear operator 
$S:L^{2}(M)\to C^{\infty}_{0}(\mathbb{R}\times M)$
such that $g=Su_{0}$ yields 
the null-controllability of the heat equation on $M$ 
in time $T$ by interior controls on $\Omega$.

The most striking feature of this result is that
we may control the heat in arbitrarily small time
whatever geometry the control region has.
In this paper we address the following question:
{\em How does the geometry of the control region 
influence the cost of controlling the heat to zero 
in small time ?}

\medskip

Now, we shall formulate this question more precisely
and give references.
\begin{defn}
\label{defin:cost}
For all control time $T$ and all control region $\Omega$,
the {\em null-controllability cost} 
for the heat equation on $M$ is the best constant, 
denoted $C_{T,\Omega}$, in the estimate:
$$
\|g\|_{L^{2}(\mathbb{R}\times M)}\leq C_{T,\Omega}\|u_{0}\|_{L^{2}(M)}
$$
for all initial data $u_{0}$ and control $g$ 
solving the null-controllability problem described above.
\end{defn}

By duality (cf.~\cite{DR77}), 
$C_{T,\Omega}$ is also the best constant in the observation inequality
for the homogeneous heat semigroup $t\mapsto e^{t\Delta}$:
$$
\forall u_{0}\in L^{2}(M),\quad
\|e^{T\Delta}u_{0}\|_{L^{2}(M)}
\leq C_{T,\Omega}
\|e^{t\Delta}u_{0}\|_{L^{2}((0,T)\times \Omega)} \ .
$$
Lebeau and Robbiano's result implies the finiteness of 
the null-controllability cost for the heat equation on $M$
for any control time and any control region.
{\`E}manuilov extended this result to more general parabolic operators 
in \cite{Ima95} using global Carleman estimates with singular weights. %cite Tataru?
When $(M,g)$ is an open set in Euclidean space, 
this method was used by Fern{\'a}ndez-Cara and Zuazua
in \cite{F-CZ00} to obtain the optimal time dependence 
of the null-controllability cost for small time, 
i.e.: %$C_{T,\Omega}=C\exp{C_{\Omega}/T}$
\begin{equation}
\label{eq:FCZ}
0<\sup_{\adh{B}_{\rho} \subset M\setminus \adh{\Omega}}\rho^{2}/4
\leq \liminf_{T\to 0} T \ln C_{T,\Omega} 
\leq \limsup_{T\to 0} T \ln C_{T,\Omega} 
< +\infty
\end{equation}
where the supremum is taken over 
balls $B_{\rho}$ of radius $\rho$.
The lower bound is stated in section 4.1 of \cite{Zua01}
%(with a misprint), 
and it is based on the construction of a 
``very singular solution of the heat equation in $(0,+\infty)\times\R^{n}$''
used in the proof of Theorem 6.2 in \cite{F-CZ00}.
Note that the method used in theorem 1 of \cite{LR95}
seems to fall short of the optimal time dependence.
Actually, using the improved version of proposition~1 in \cite{LR95} 
presented as proposition~2 in \cite{LZ98}, 
we have only been able to prove that
$\limsup_{T\to 0} T^{\gamma} \ln C_{T,\Omega}$ 
is finite for all $\gamma>1$.

Indeed Seidman had already asked how violent fast controls are,
and his first answer concerned heat null-controllability
from a boundary region $\Gamma\subset \partial M$.
In \cite{Sei84}, 
under the condition that the wave equation on $M$ is 
exactly controllable by controls in $\Gamma$ in time $L$,
he computes an explicit positive value $\beta$
such that $\limsup_{T\to 0} T \ln C_{T,\Gamma}\leq \beta L^{2}$
(we give more explanations on this geometric upper bound 
in section~\ref{sec:res} after theorem~\ref{theo:ub}).
The positivity of $\liminf_{T\to 0} T \ln C_{T,\Gamma}$ 
when $M$ is an interval was subsequently proved by G{\"u}ichal  
in \cite{Gui85}, 
ensuring the optimality of Seidman's result 
with respect to the time dependence. %in dimension $n=1$. 
Later, Seidman also addressed finite dimensional linear systems %in~\cite{Sei88}
as well as the Schr{\"o}dinger and plate equations 
(cf.~the companion paper \cite{LMschrocost}
for more details and references).

%%%%%%%%%%%%%%%%%%%%%%%%%%%%%%%%%%%%%%%%%%%%%%%%%%%%%%%%%%%%
%% RESULTS
%%%%%%%%%%%%%%%%%%%%%%%%%%%%%%%%%%%%%%%%%%%%%%%%%%%%%%%%%%%%

\section{The results}
\label{sec:res}

\subsection{Lower bound}
Our first result, proved in section~\ref{sec:lb}, 
generalizes and improves 
on the geometric lower bound of Fern{\'a}ndez-Cara and Zuazua:
\begin{thm}
\label{theo:lb}
The null-controllability cost of the heat equation for small time
(cf.\ definition~\ref{defin:cost}) 
satisfies the following geometric lower bound:
\begin{equation}
\label{eq:lb}
\liminf_{T\to 0} T \ln C_{T,\Omega} \geq \sup_{y\in M}\dist(y,\adh{\Omega})^{2}/4
\end{equation}
\end{thm}

As put in \cite{Zua01},
such a lower bound follows from 
the construction of a ``very singular solution of the heat equation''.
Our construction underscores that 
only a large but finite number of modes is needed. 
For a short control time $T>0$,
we consider a Dirac mass as far from $\Omega$ as possible, 
we smooth it out by applying 
the homogeneous heat semigroup for a very short time
($\eps T$ with small $\eps$)
and truncating very large frequencies 
(larger than $(\eps T)^{-1}$),
and finally we take it as initial data in (\ref{eqHeat}).
The proof relies on Varadhan's formula 
for the heat kernel in small time (cf.~\cite{Var67}),
which requires very low smoothness assumptions 
as proved in~\cite{Nor97}.

We believe that there is no solution of the heat equation
which is more singular than the heat kernel
and therefore conjecture that this lower bound is also an upper bound, 
i.e. $\displaystyle 
\lim_{T\to 0} T \ln C_{T,\Omega} = \sup_{y\in M}\dist(y,\adh{\Omega})^{2}/4$.

\subsection{The segment controlled at one end}
Our second result, proved in section~\ref{sec:1d}, 
concerns the most simple heat null-controllability problem:
the heat equation on a segment controlled at one end
through a Dirichlet condition.
It is an upper bound of the same type as 
the lower bound in theorem~\ref{theo:lb},
except that the quite natural rate $1/4$
is replaced by the technical rate
(resulting from the complex multiplier lemma~\ref{lemM}):
\begin{equation}
\label{eqalphastar}
\alpha_{*}=2\left(\frac{36}{37}\right)^{2}
%% 2\min\left\{  \left(\frac{36}{37}\right)^{2},
%% \left(1+\sum_{k\in \mathbb{N}^{*}} 
%% \frac{1}{2k(4k-1)}
%% \frac{\zeta(2k)}{\pi^{2k}}
%% \right)^{-1}
%% \right\}
<2
\ .
\end{equation}
%where %$\zeta$ denotes the zeta function:
%$\zeta(s)=\sum_{n\in \mathbb{N}^{*}}n^{-s}$.
\begin{thm}
\label{theo:1d}
For any $\alpha > \alpha_{*}$ defined by $(\ref{eqalphastar})$,
there exists $C>0$ such that, 
for $B=1$ or $B=\partial_{s}$,
for all $L>0$, $T\in\, ]0,\inf(\pi,L)^{2}]$
and $u_{0}\in L^{2}(0,L)$,
there is a $g\in L^{2}(0,T)$
such that the solution $u\in C^{0}([0,\infty),L^{2}(0,L))$
of the following heat equation on $[0,L]$
controlled by $g$ from one end: 
\begin{equation*} 
\label{eqHeat1d}
\partial_{t}u - \partial_{s}^{2} u=0
\quad {\rm in}\ ]0,T[\times ]0,L[\, , \quad
\left( Bu \right)\res{s=0} = 0 \, ,\quad 
u\res{s=L} = g \, ,\quad
u\res{t=0} = u_{0} \, ,
\end{equation*}
satisfies $u=0$ at $t=T$ and 
$\displaystyle
\|g\|_{L^{2}(0,T)}\leq Ce^{\alpha L^{2}/T }\|u_{0}\|_{L^{2}(0,L)} \ .
$
\end{thm}
Theorem~3.1 in \cite{Sei84}
yields this theorem for $\alpha_{*}=4\beta_{*}$
with $\beta_{*}\approx 42.86$.
This result of Seidman can be improved to $\alpha_{*}=8\beta_{*}$ 
with $\beta_{*}\approx 4.17$
using his theorem~1 in \cite{Sei86}.
%Following \cite{SAI00}, 
%which applies to much more general spectral sequences,
%yields this theorem for some $\alpha_{*}$ greater than $4\times 12$
%(cf. the remarks after lemmas~\ref{lemF} and~\ref{lemM}).
The value $\alpha_{*}$ defined by $(\ref{eqalphastar})$ 
in theorem~\ref{theo:1d}
is the best we obtained yet  
following the well trodden path 
of the harmonic analysis of this problem
(cf.~\cite{Rus78} and \cite{SAI00} for seminal and recent references).
As explained at the end of the previous subsection, 
we conjecture that $\alpha_{*}=1/4$ 
is the optimal rate.
The proof of theorem~\ref{theo:lb} also applies here, 
so that  theorem~\ref{theo:1d} does not hold with $\alpha_{*}<1/4$. 
This theorem is valid 
for more general linear parabolic equations and boundary conditions
as formulated in theorem~\ref{theo:1dParabolic}.

\subsection{Upper bound under the geodesics condition}
Our third result gives a good reason to strive for the best rate 
$\alpha_{*}$ in theorem~\ref{theo:1d}.
In section~\ref{sec:transmut}, we prove that the upper bound 
for the null-controllability cost 
of the heat equation on a segment controlled at one end
---~the particular case in which the computation are the most explicit~---
is also an upper bound for the multidimensional case 
of equation (\ref{eqHeat})
under the following {\em geodesics condition} on the control region:
every generalized geodesic
in $\adh{M}$
intersects $\Omega$.

In this context, the {\em generalized geodesics} 
are continuous trajectories $t\mapsto x(t)$ in $\adh{M}$ 
which follow geodesic curves at unit speed in $M$
(so that on these intervals 
$t\mapsto \dot{x}(t)$ is continuous);
if they hit $\partial M$ transversely at time $t_{0}$,
then they reflect as light rays or billiard balls
(and $t\mapsto \dot{x}(t)$ 
is discontinuous at $t_{0}$); 
if they hit $\partial M$ tangentially 
then either there exists a geodesic in $M$
which continues $t\mapsto (x(t),\dot{x}(t))$
continuously 
and they branch onto it,
or there is no such geodesic curve in $M$
and 
then they glide at unit speed 
along the geodesic of $\partial M$
which continues $t\mapsto (x(t),\dot{x}(t))$
continuously until they may branch onto  a geodesic in $M$.
For this result and whenever generalized geodesics are mentionned, 
we make the additional assumptions
that they can be uniquely continued 
at the boundary $\partial M$
(as in~\cite{BLR92}, to ensure this, we may assume either that 
$\partial M$ has no contacts of infinite order with its tangents,
or that $g$ and $\partial M$ are real analytic),
and that $\Omega$ is open.

\begin{thm}
\label{theo:ub}
Let $L_{\Omega}$ be the length of the longest generalized geodesic in $\adh{M}$
which does not intersect $\Omega$.
If theorem~\ref{theo:1d} holds for some rate $\alpha_{*}$
then the null-controllability cost of the heat equation for small time 
(cf.\ definition~\ref{defin:cost}) 
satisfies the following geometric upper bound:
\begin{equation}
\label{eq:ub}
\limsup_{T\to 0} T \ln C_{T,\Omega} \leq \alpha_{*}
L_{\Omega}^{2}
\end{equation}
\end{thm}

When comparing this result to the lower bound in theorem~\ref{theo:lb},
one should bear in mind that $L_{\Omega}$
is always greater than $2\sup_{y\in M}\dist(y,\adh{\Omega})$
(because the length of a generalized geodesic through $y$
which does not intersect $\Omega$ 
is always greater than $2\dist(y,\adh{\Omega})$)
%(take $y$ half way on the longest generalized geodesic in $\adh{M}$)
and can be infinitely so.
For instance,
on the sphere $M=S^{n}$,
if $\Omega$ is the complementary set 
of a tube of radius $\eps$ around the equator, 
then 
$\sup_{y\in M}\dist(y,\adh{\Omega})=\eps$
and $L_{\Omega}=\infty$.
If $\Omega$ is increased by 
a tube slice of small thickness $\delta$,
then the first length is unchanged
while the second length becomes 
greater than the length of the equator of $M$ minus $\delta$,
so that $L_{\Omega}$ is finite 
yet much greater than $\sup_{y\in M}\dist(y,\adh{\Omega})$
as $\eps\to 0$.

Moreover, 
as recalled in section~\ref{sec:pb},
this geodesics condition 
is by no means necessary for 
the null-controllability of the heat equation. 
It is more relevant to the wave equation on $M$,
for which it is a sharp sufficient condition 
for exact controllability in time $T$ by interior controls on $\Omega$
as proved in~\cite{BLR92}
(cf. theorem~\ref{theo:BLR} for the precise statement).
It was later proved in~\cite{BG97} 
that this condition is also necessary when 
the characteristic function of $]0,T[\times \Omega$ 
is replaced by a smooth function $\theta$
such that $\{\theta(t,x)\neq 0\}=]0,T[\times \Omega$. 
%% \footnote{
%% In particular, if $L_{\Omega}=+\infty$
%% (e.g. if there is a periodic generalized geodesic
%% which does not intersect $\Omega$), 
%% then the wave equation on $M$ 
%% is not exactly controllable 
%% by interior controls on $\Omega$
%% for any time $T$.
%% }.

In fact we use the exact controllability of the wave equation
to prove our result on the null-controllability of the heat equation. 
This strategy was already applied by Russell in 1973, %\cite{Rus73}  
but he used a complex analysis detour (cf.~\cite{Rus78}).
In \cite{Sei84}, Seidman applied Russell's method
to obtain an upper bound %(already mentionned in section~\ref{sec:pb})
which, taking~\cite{BLR92} into account, 
corresponds to theorem~\ref{theo:ub} 
with $\alpha_{*}=\beta_{*}\approx 42.86$.
%for $\alpha_{*}=\beta_{*}$ with $\beta_{*}\approx 42.86$.
Theorem~\ref{theo:ub} improves Seidman's result
beyond this slight improvement of the rate $\alpha_{*}$ 
insofar as the complex analysis multiplier method he uses 
does not necessarily allow to reach 
the optimal $\alpha_{*}$ in theorem~\ref{theo:1d}.

The {\em control transmutation method}
(cf.~\cite{Her75} for a survey on transmutations in other contexts)
% (about transmutation methods)
introduced in section~\ref{sec:transmut}  
relates the null-controllability of the heat equation
to the exact controllability of the wave equation
in a direct way 
(as opposed to Russell's indirect complex analysis method).
It is well-known that 
the geometry of small time asymptotics 
for the homogeneous heat semigroup $t\mapsto e^{t\Delta}$ on $L^{2}(M)$
can be understood from the even homogeneous wave group $t\mapsto W(t)$
(i.e. the group defined %for all $w_{0}\in L^{2}(M)$
by $W(t)w_{0}=w(t)$ where $w$ solves equation (\ref{eqWave})
with $f=0$ and Cauchy data $(w,\partial_{t}w)=(w_{0},0)$ at $t=0$)
through Kannai's formula 
(cf.~\cite{Kan77}, \cite{CGT82}, and section 6.2 in the book~\cite{Tay96})~:
\begin{equation} 
\label{eqKannai}
e^{t\Delta}=\frac{1}{\sqrt{4\pi t}}
\int_{-\infty}^{\infty} e^{-s^{2}/(4t)}W(s)\, ds 
\ .
\end{equation}
Our main idea is to replace the fundamental solution 
of the heat equation on the line 
$e^{-s^{2}/(4t)}/\sqrt{4\pi t}$
appearing in Kannai's formula 
by some  {\em fundamental controlled solution} 
of the heat equation on the segment $[-L,L]$ 
controlled at both ends.
We use the one dimensional theorem~\ref{theo:1d} to 
construct this fundamental controlled solution
in section~\ref{sec:transmut}.

\subsection{Open problems}
We shall now survey some questions raised by the results 
we have presented which we have been unable to answer yet.   

To improve the rate $\alpha_{*}$ in theorem~\ref{theo:1d}
by a complex analysis method,
one could use the first method in \cite{FR71},
i.e. compute the null-controllability 
cost on the half-line $[0,+\infty)$ explicitly 
by Vandermonde determinants 
and prove a quantitative version of Schwartz's theorem in \cite{Sch43},
i.e. estimate with respect to $L$ the best constant $c_{L}$ 
in the following statement~: 
every $u$ in the closed linear hull in $L^{2}(0,+\infty)$
of the real exponential sums $t\mapsto e^{-k^{2}t}$ ($k\in \mathbb{N}^{*}$) 
satisfies
$\displaystyle
\| u \|_{L^{2}(0,+\infty)}\leq c_{L} \| u \|_{L^{2}(0,L)}
$.

Theorem~\ref{theo:ub} opens new tracks 
to improve the upper bound for 
the null-controllability cost of $(\ref{eqHeat})$ 
under the geodesics condition
by methods which are not complex analytical.
To improve the rate $\alpha_{*}$ in theorem~\ref{theo:1d}
(or in the multidimensional case of equation (\ref{eqHeat})
when $\Omega$ and $M$ are star-shaped with respect to the same point)
one could adapt the variational techniques 
(e.g. the log convexity method)
or the Carleman's inequalities 
devised to prove unique continuation theorems.
 
In the general case (without the geodesics condition),
one could try to adapt the null-controllability proofs
which use Carleman inequalities with phases $\phi$
to obtain an upper bound similar to the lower bound in theorem~\ref{theo:lb}
in terms of the following distance function $d$~:
$d(x,y)=\sup\{\phi(y)-\phi(x)\}$,
for all $x$ and $y$ in $M$, 
where the supremum is taken over all Lipschitz functions 
$\phi:M\to \mathbb{R}$ with $|\nabla\phi|\leq 1$ almost everywhere.
There is a more geometric characterization of %this distance function 
$d$
in terms of path of least action (cf.~section~2 of \cite{Nor97}).

%%%%%%%%%%%%%%%%%%%%%%%%%%%%%%%%%%%%%%%%%%%%%%%%%%%%%%%%%%%%
%% LOWER BOUND
%%%%%%%%%%%%%%%%%%%%%%%%%%%%%%%%%%%%%%%%%%%%%%%%%%%%%%%%%%%%

\section{Lower bound}
\label{sec:lb}

The purpose of this section is to prove theorem~\ref{theo:lb}.

As in section~\ref{sec:pb}, 
let $\Omega$ be an open set in the 
$n$-dimensional Riemannian manifold $M$
such that $\adh{\Omega}\subset M$.
Let $(\omega_{j})_{j\in \mathbb{N}^{*}}$
be a nondecreasing sequence of nonnegative real numbers 
and $(e_{j})_{j\in \mathbb{N}^{*}}$
be an orthonormal basis of $L^{2}(M)$
such that $e_{j}$ is an eigenvector of $-\Delta$
with eigenvalue $\omega_{j}^{2}$.
The heat kernel $k$ can be defined 
for all $t>0$ and $(x,y)\in \adh{M}^{2}$ by 
$k(t,x,y)=\sum_{j}\exp(-t\omega_{j}^{2})e_{j}(y)e_{j}(x)$.
Our main ingredient is Varadhan's formula which says that 
(cf. th.~1.1 in \cite{Nor97} for example):
\begin{equation}
  \label{eqVaradhan}
  \lim_{t\to 0} t\ln k(t,x,y) = -d(x,y)^{2}/4 \quad
\mbox{ uniformly on compact sets of } \adh{M}^{2}
\ .
\end{equation}
We shall also use Weyl's asymptotics for eigenvalues:
\begin{equation}
  \label{eqWeyl}
  \exists W>0,\  
\#\{j\in \mathbb{N}^{*}\, | \, \omega_{j}\leq \omega \} \leq W \omega^{n}  
\end{equation}
and the following consequence of Sobolev's embedding theorem:
\begin{equation}
  \label{eqSobolev}
  \exists E>0,\,  \forall  j\in \mathbb{N}^{*},\
\|e_{j}\|_{L^{\infty}}\leq  E \omega_{j}^{n/2}
\end{equation}
(cf. section 17.5 in \cite{HIII} for example).
The unique continuation property for elliptic operators implies that 
$Y=\{ y\in M\setminus \adh{\Omega} \, |\, e_{1}(y)\neq 0 \}$
is an open dense set in $M\setminus\adh{\Omega}$,
so that the supremun in theorem~\ref{theo:lb}
can be taken over $y\in Y$ instead of $y\in M$.

Let $y\in Y$ and $\alpha < d(y,\adh{\Omega})^{2}/4$
be fixed from now on.
To prove theorem~\ref{theo:lb} we shall find $A>0$ 
and, for all $T\in ]0,1]$ small enough,
some data $u_{0}^{T}\in L^{2}(M)$ such that 
$\|e^{t\Delta}u_{0}^{T}\|_{L^{2}((0,T)\times \Omega)}
\leq 
Ae^{-\alpha/T}
\|e^{T\Delta}u_{0}^{T}\|_{L^{2}(M)}
$.
To give further insight into the problem, 
we shall construct each $u_{0}^{T}$ 
as a linear combination of
a finite number of modes $e_{j}$ only.

Let $\beta$ be a real number such that 
$\alpha<\beta<d(y,\adh{\Omega})^{2}/4$.
Since $\adh{\Omega}\times \{y\}$ is compact in $\adh{M}^{2}$,
Varadhan's formula (\ref{eqVaradhan}) yields real numbers 
$B>0$ and $\Tb\in ]0,1]$ such that 
\begin{equation}
  \label{eqk}
\forall t\in ]0,\Tb],\, 
\forall x\in \adh{\Omega}, \
|k(t,x,y)|\leq Be^{-\beta/t}
\ .
\end{equation}
Let $\eps\in ]0,1]$ small enough as specified later. 
For all $T\in ]0, \Tb/(1+\eps)]$
consider the data 
$u_{0}^{T}(x)=\sum_{\omega_{j}\leq (\eps  T)^{-1}}
\exp(-\eps T\omega_{j}^{2})e_{j}(y)e_{j}(x)$.
To estimate the corresponding solution 
$$
u^{T}(t,x)=\left( e^{t\Delta}u_{0}^{T} \right)(x)
=\sum_{\omega_{j}\leq (\eps  T)^{-1}}
\exp(-(\eps T+t)\omega_{j}^{2})e_{j}(y)e_{j}(x)
\ ,
$$
we compare it with $k(\eps T+t,x,y)$.
Using that the heat semigroup is a contraction on $L^{2}(M)$,
Parseval's identity and (\ref{eqSobolev}), 
we obtain
\begin{eqnarray*}
\sup_{t\in ]0,T]} \| k(\eps T+t,x,y)-u^{T}(t,x)\|_{L^{2}(M)}
\leq  \| k(\eps T,x,y)-u^{T}_{0}(x)\|_{L^{2}(M)} \\
=\sum_{\omega_{j}> (\eps  T)^{-1}}
|e^{-\eps T\omega_{j}^{2}}e_{j}(y)|^{2}
\leq E \sum_{\omega_{j}\geq (\eps  T)^{-1}}
e^{-\omega_{j}}\omega_{j}^{n}
\leq %2^{n}
E'\sum_{\omega_{j}\geq (\eps  T)^{-1}}e^{-\omega_{j}/2}
\ ,
\end{eqnarray*}
for some $E'>0$.
But, Weyl's law (\ref{eqWeyl}) yields, 
for $c\geq c_{0}>0$ and $\gamma\geq \gamma_{0}>0$,
\begin{eqnarray*}
\sum_{\omega_{j}\geq c}e^{-\gamma \omega_{j}}
=\sum_{k\in \mathbb{N}^{*}} 
\sum_{kc\leq \omega_{j} <(k+1)c}e^{-\gamma \omega_{j}}
\leq W \sum_{k\in \mathbb{N}^{*}} 
\left( (k+1)c \right)^{n} e^{-kc\gamma } \\
%\leq  \left(\frac{4}{\gamma}\right)^{n}W' 
\leq  W_{\gamma_{0}} 
\sum_{k\in \mathbb{N}^{*}} e^{-kc\gamma}e^{(k+1)c\gamma/4}
= W_{\gamma_{0}} e^{-c\gamma/2}
\sum_{k\in \mathbb{N}} e^{-3kc\gamma/4} 
\leq W_{c_{0},\gamma_{0}} e^{-c\gamma/2}
\end{eqnarray*}
where $W_{\gamma_{0}}$ and $W_{c_{0},\gamma_{0}}$
are positive real numbers which depend on their indices 
but not on $c$ and $\gamma$.
Hence, with $c=(\eps  T)^{-1}>1=c_{0}$ and $\gamma=\gamma_{0}=1/2$, we obtain: 
\begin{equation*}
 \exists B'>0,\ 
 \forall t\in ]0,T]\ 
 \| k(\eps T+t,x,y)-u^{T}(t,x)\|_{L^{2}(M)}
\leq B' e^{-1/(4\eps T)}
\end{equation*}
Together with the estimate on $k(\eps T+t,x,y)$
which follows from (\ref{eqk}),  
this estimate yields by the triangle inequality,
choosing $\eps<1/(4\beta)$
and setting $B''=|\Omega|^{1/2}B + B'$,
$$
\| u^{T} \|_{L^{2}((0,T)\times \Omega)}
\leq (T|\Omega|)^{1/2}Be^{-\beta/((1+\eps)T)}
+ T^{1/2}B' e^{-1/(4\eps T)}
\leq B''e^{-\beta/((1+\eps)T)}
\ .
$$
But using Parseval's identity and $y\in Y$, we have
$$
\|e^{T\Delta}u_{0}^{T}\|_{L^{2}(M)}
%\|u^{T}(T,x)\|_{L^{2}(M)}
=\left( 
\sum_{\omega_{j}\leq (\eps  T)^{-1}}
|e^{-(1+\eps) T\omega_{j}^{2}}e_{j}(y)|^{2}
\right)^{1/2}
\geq e^{-2\omega_{1}^{2}} |e_{1}(y)|
>0
\ .
$$
Hence, 
choosing $\eps$ small enough so that $\alpha < \beta/(1+\eps)$
and setting \linebreak $A=e^{-2\omega_{1}^{2}} |e_{1}(y)|B''$,
we have  
$$
\forall T\in ]0, \Tb/(1+\eps)],\ 
\|u^{T}\|_{L^{2}((0,T)\times \Omega)}
\leq Ae^{-\alpha/T}\|e^{T\Delta}u_{0}^{T}\|_{L^{2}(M)}
\ .
$$
Since $A$ does not depend on $T$, 
this ends the proof of theorem~\ref{theo:lb}.

%%%%%%%%%%%%%%%%%%%%%%%%%%%%%%%%%%%%%%%%%%%%%%%%%%%%%%%%%%%%
%% ONE DIMENSIONAL PB
%%%%%%%%%%%%%%%%%%%%%%%%%%%%%%%%%%%%%%%%%%%%%%%%%%%%%%%%%%%%

\section{The segment controlled at one end}
\label{sec:1d}
%{Harmonic analysis of the segment controlled at one end}
%{One dimensional problem}

In this section we prove theorem~\ref{theo:1d} 
for a more general linear parabolic equation
on a segment controlled at one end
(in particular, 
it proves that theorem~\ref{theo:1d}
is true for the heat equation 
on a segment with any Riemannian metric).
We follow \cite{FR71} quite closely. 

For a positive a control time $T$,
we consider the following mixed Dirichlet-Cauchy problem
on the space segment $[0,X]$:
\begin{eqnarray} 
\label{eqParab1}
&&
\partial_{t}u= \partial_{x}\left(p(x)\partial_{x}u\right) +q(x)u
\quad {\rm for}\ (t,x)\in ]0,T[\times ]0,X[\, , %\quad 
\\ 
\label{eqParab2}
&&
\left(a_{0}+b_{0}\partial_{x}\right)u\res{x=0} = 0 \, ,\quad
\left(a_{1}+b_{1}\partial_{x}\right)u\res{x=X} = g \, ,\quad
u\res{t=0} = u_{0} \ ,
\\ 
\label{eqParabAssum}
&&
a_{0}^{2}+b_{0}^{2}= a_{1}^{2}+b_{1}^{2}=1 \, ,\quad
0<p\in C^{2}([0,X]) % \ \mbox{is positive} 
\, ,\quad
q\in C^{0}([0,X])
\ . 
\end{eqnarray}
%We assume $a_{0}^{2}+b_{0}^{2}= a_{1}^{2}+b_{1}^{2}=1$,
%$p\in C^{2}([0,X])$ is positive and $q\in C^{0}([0,X])$.
With assumptions $(\ref{eqParabAssum})$, 
the operator $A$ on $L_{2}(0,X)$ with domain $D(A)$  
defined by 
\begin{eqnarray*}
(Au)(x) &=& \partial_{x}\left(p(x)\partial_{x}u(x)\right) +q(x)u(x)
\\
D(A) &=& H^{2}(0,X)\cap \{
\left(a_{0}+b_{0}\partial_{x}\right)u\res{x=0} = 
\left(a_{1}+b_{1}\partial_{x}\right)u\res{x=X} = 0
\}
\end{eqnarray*}
is self-adjoint
and has a sequence $\{-\lambda_{n}\}_{n\in \mathbb{N}^{*}}$
of increasing eigenvalues 
and an orthonormal Hilbert basis $\{e_{n}\}_{n\in \mathbb{N}^{*}}$ 
in $L_{2}(0,X)$ of corresponding eigenfunctions, i.e. : 
\begin{equation*}
\forall n\in \mathbb{N}^{*}, \quad -Ae_{n}  = \lambda_{n}e_{n} 
\quad \mbox{and} \quad 
\lambda_{n}<\lambda_{n+1} \ .
\end{equation*}
Moreover, $(\ref{eqParabAssum})$ ensures the following eigenvalues asymptotics
%there is a real constant $\nu$ such that 
(cf. \cite{FR71}):
\begin{equation}
\label{eqParabEig}
\exists \nu\in \mathbb{R},\
\lambda_{n}=\frac{\pi^{2}}{L}\left(n+\nu\right)^{2} + O(1) 
\ \mbox{as} \  n\to \infty \ ,
\quad \mbox{where} \quad L=\int_{0}^{X}\sqrt{p(x)}\, dx
\ .
\end{equation}

\begin{thm}
\label{theo:1dParabolic}
For any $\alpha > \alpha_{*}$ defined by $(\ref{eqalphastar})$,
there exists $C>0$ such that, 
for any coefficients $(\ref{eqParabAssum})$, 
%such that $L\geq \sqrt{T}$,
%for all $T\in\, ]0,1[$, $L\geq \sqrt{T}$,
for all $T\in\, ]0,\inf(\pi,L)^{2}]$
and $u_{0}\in L^{2}(0,X)$
there is a control $g\in L^{2}(0,T)$
such that the solution $u\in C^{0}([0,\infty),L^{2}(0,X))$
of $(\ref{eqParab1})$ and $(\ref{eqParab2})$ 
satisfies $u=0$ at $t=T$ and 
$\displaystyle
\|g\|_{L^{2}(0,T)}\leq Ce^{\alpha L^{2}/T }\|u_{0}\|_{L^{2}(0,X)} \ .
$
%with $L=\int_{0}^{X}\sqrt{p(x)}\, dx$.
\end{thm}

As in \cite{FR71}, the proof applies to the slightly more general  
eigenvalue asymptotics 
$\lambda_{n}=\frac{\pi^{2}}{L}\left(n+\nu\right) + o(n)$.
We divide the proof of this theorem in three steps.

\subsection{Reduction to positive eigenvalues, to a segment of $p$-length $L=\pi$, 
and to the control window $]-T/2,T/2[$}

As a first step, 
we reduce the problem to the case 
$\lambda_{1}>0$ 
by the multiplier $t\mapsto \exp(\lambda t)$,
to the case $L=\pi$ 
by the time rescaling $t\mapsto \sigma t$
with $\sigma=(\pi/L)^{2}$,
and to the time interval $[-T/2,T/2]$
by the time translation $t\mapsto t-T/2$.

The function $u$ satisfies $\partial_{t}u=Au$
and $\left(a_{1}+b_{1}\partial_{x}\right)u\res{x=X} = g$
if and only if $\ut(t,x)=\exp(\lambda t)u(t,x)$
satisfies $\partial_{t}\ut=\At \ut$
and $\left(a_{1}+b_{1}\partial_{x}\right)\ut\res{x=X} = \gt$
with $\At=A+\lambda$ and $\gt(t)=\exp(\lambda t)g(t)$.
For any $\lambda>-\lambda_{1}$, 
the lowest eigenvalue of 
$\At\geq \lambda_{1}+\lambda>0$
is positive.
In $\At$, $q$ is changed into $q+\lambda$
and $p$ is unchanged so that $L$ is unchanged.
Moreover 
$\|g\|_{L^{2}(0,T)}\leq \exp(\lambda T/2)\|\gt\|_{L^{2}(0,T)}$
so that 
$
\|\gt\|_{L^{2}(0,T)}\leq \Ct e^{\alpha L^{2}/T }\|u_{0}\|_{L^{2}(0,X)} 
$
implies the estimate in theorem~\ref{theo:1dParabolic} 
with \linebreak $C=\Ct \exp(\lambda \pi/2)$.
This proves the reduction to positive eigenvalues.

We now prove the second reduction.
Assume the theorem is true when $L$ takes the value $\Lt=\pi$.
Given $L>0$ and $T\in\, ]0,\inf(\pi,L)^{2}]$
we set $\Tt=\sigma^{2}T\in ]0,\Lt^{2}]$ and $\At=\sigma^{2}A$,
where $\sigma=(\pi/L)^{2}$.
By applying the theorem to $\At$ on $]0,\Tt[$,
we obtain  
$%\displaystyle
\|\gt\|_{L^{2}(0,\Tt)}\leq \Ct e^{\at \Lt^{2}/\Tt }\|u_{0}\|_{L^{2}(0,X)} 
$.
The function $g(t)=\gt(\sigma t)$
is a control for the solution $u(t,x)=\ut(\sigma t,x)$ 
of $\partial_{t} u = A u $ on $]0,T[$
at the cost 
$\|g\|_{L^{2}(0,T)}=\|\gt\|_{L^{2}(0,T)}L/\pi$.
Since $T\leq \pi^{2}$ implies $L/\pi\leq (L^{2}/T)^{1/2}$,
for all $\alpha>\at$ there is a $C$ such that 
for all $L>0$ and $T\in\, ]0,\inf(\pi,L)^{2}]$: 
$\Ct e^{\at \Lt^{2}/\Tt }L/\pi\leq Ce^{\alpha L^{2}/T }$.
Therefore $g$ satisfies the estimate in theorem~\ref{theo:1dParabolic}.

These two reductions allow us to assume from now on 
$\lambda_{1}>0$ and $L=\pi$.
Making a weaker assumption on the remainder term in $(\ref{eqParabEig})$,
we shall only use the following spectral assumption:
\begin{equation}
\label{eqSpec}
\forall n\in \mathbb{N}^{*},\  
0<\lambda_{n}<\lambda_{n+1} 
\quad \mbox{and} \quad 
\exists \nu\in \mathbb{R},\ 
\lambda_{n}=\left(n+\nu\right)^{2} + o(n) 
\ \mbox{as} \ n\to \infty \ .
\end{equation}

It is obvious that theorem~\ref{theo:1dParabolic}
is invariant by time translations
and we shall prove it for the control window $]-T/2,T/2[$
instead of $]0,T[$.

\subsection{Spectral reduction to a problem in complex analysis}

In this second step, 
we recall that the control $g$ in this theorem can be obtained 
as a series expansion into a Riesz sequence 
$\{g_{n}\}_{n\in \mathbb{N}^{*}}$ in $L^{2}(-T/2,T/2)$
which is bi-orthogonal to the sequence 
$\{\exp(-\lambda_{n}t)\}_{n\in \mathbb{N}^{*}}$.
We also recall 
how the Paley-Wiener theorem reduces 
the construction of such biorthogonal functions
to the construction of entire functions 
with zeros and growth conditions
(this well-known method in complex analysis 
is the second method in \cite{FR71}
called the Fourier transform method there). 
Our estimate on the control cost $\|g\|_{L^{2}(-T/2,T/2)}$ 
relies on a good estimate of $\|g_{n}\|_{L^{2}(-T/2,T/2)}$ 
as $T$ tends to zero.
This additional difficulty was first taken care of 
by Seidman in \cite{Sei86} for $\lambda_{n}=in^{2}$
and it was recently overcome for more general sequences in \cite{SAI00}.
Our contribution is a slight improvement on 
the estimates of Seidman and his collaborators
in our less general setting. 

In terms of the coordinates
$c=(c_{j})_{j\in \mathbb{N}^{*}}$ of $u_{0}$ 
in the Hilbert basis $(e_{j})_{j\in \mathbb{N}^{*}}$,
the controllability problem in theorem~\ref{theo:1dParabolic}
is equivalent to the following moment problem 
(by straightforward integration by parts, cf. \cite{FR71}):
$$%\displaystyle
\int_{-T/2}^{T/2}e^{-\ldn (T/2-t)}\gamma_{n} g(t) \, dt=- e^{-\ldn T}c_{n}
\ ,$$
where $\gamma_{n}=e_{n}(X)p(X)/b_{1}$ if $b_{1}\neq 0$
and  $\gamma_{n}=-e_{n}'(X)p(X)/a_{1}$ if $b_{1}= 0$.
In both cases, the asymptotic expansion of $e_{n}$ 
yields %$0\neq\gamma_{n}\sim c\sqrt{\ldn}$ (cf. \cite{FR71}), so 
that $(|\gamma_{n}|)$ is bounded from below 
by some positive constant $\gamma$.
If $\{g_{n}\}_{n\in \mathbb{N}^{*}}$ in $L^{2}(-T/2,T/2)$
is a sequence which is bi-orthogonal to the sequence 
$\{\exp(-\lambda_{n}t)\}_{n\in \mathbb{N}^{*}}$, i.e.
%(i.e. $\int_{0}^{T}g_{n}(t)e^{-\ldk t} \, dt=\delta_{n,k}$),
\begin{equation}
\label{eqbiorth}
\int_{-T/2}^{T/2}g_{n}(t)e^{-\ldn t} \, dt=1 
\quad \mbox{ and  } \quad
\forall k\in \mathbb{N}^{*} ,\, k\neq n ,\,  
\int_{-T/2}^{T/2}g_{n}(t)e^{-\ldk t} \, dt=0  , 
\end{equation}
then 
$\displaystyle
g(t)=-\sum_{n=1}^{\infty} \frac{c_{n}}{\gamma_{n}}e^{-\ldn T/2} g_{n}(-t)
$
is a formal solution to this moment problem.
The following theorem in complex analysis
allows to construct a bi-orthogonal sequence
such that this series converges 
and yields a good estimate of $\|g\|_{L^{2}(-T/2,T/2)}$  
as $T$ tends to zero.
% To prove theorem~\ref{theo:1dParabolic}, 
% we need this series to converge 
% and a good estimate of its norm as $T$ tends to zero.
% Therefore we need to construct the biorthogonal sequence 
% with a good estimate of $\|g_{n}\|_{L^{2}(0,T)}$ 
% as $T$ tends to zero.
% the Paley-Wiener theorem reduces this construction 
% to the following theorem in complex analysis.
\begin{thm}
  \label{theo:G}
Let $\alpha_{*}$ be defined by $(\ref{eqalphastar})$.
Let $\{\lambda_{n}\}_{n\in \mathbb{N}^{*}}$ be a sequence of real numbers 
satisfying $(\ref{eqSpec})$.
For all $\eps>0$
%there are positive constants $b$ and $B$ such that,
there is a $C_{\eps}>0$ such that, 
for all $\tau\in ]0,1]$
and $n\in \mathbb{N}^{*}$, 
there is an entire function $G_{n}$ 
satisfying
\begin{eqnarray}
&&
  G_{n} \mbox{ is of exponential type }\tau
\mbox{, i.e. } \limsup_{r\to +\infty} r^{-1}\sup_{|z|=r}\ln |G_{n}(z)| \leq\tau,
\label{eqG1} \\
&&
G_{n}(i\ldn)=1 
\quad \mbox{ and  } \quad
\forall k\in \mathbb{N}^{*} ,\, k\neq n ,\,  
G_{n}(i\ldk)=0 , 
\label{eqG2} \\
&&
\|G_{n}\|_{L^{2}}=\left(
\int_{-\infty}^{+\infty} |G_{n}(x)|^{2} \, dx
\right)^{1/2} 
\leq  C_{\eps}e^{\eps \sqrt{\ldn}}e^{\alpha_{*}
(\pi+2\eps)^{2}/(2\tau)}
\label{eqG3} 
\end{eqnarray}
\end{thm}

According to the Paley-Wiener theorem (1934),
$(\ref{eqG1})$ implies that the function $x\mapsto G_{n}(x)$
is the unitary Fourier transform of a function $t\mapsto g_{n}(t)$
in $L^{2}(\mathbb{R})$ supported in $[-\tau,\tau]$.
With $\tau=T/2$, this yields:
\begin{equation}
\label{eqg}
G_{n}(x)=\frac{1}{\sqrt{2\pi}}\int_{-T/2}^{T/2}g_{n}(t)e^{-itx}\, dt
\quad \mbox{and} \quad 
\|g_{n}\|_{L^{2}}=\|G_{n}\|_{L^{2}} \ .
\end{equation}
Hence $(\ref{eqG2})$ implies $(\ref{eqbiorth})$   
and $(\ref{eqG3})$ implies that the series defining $g$
converges with:
$$
\|g\|_{L^{2}}\leq 
\sum_{n=1}^{\infty} \left| \frac{c_{n}}{\gamma_{n}} \right| 
e^{-\ldn T/2} 
\|g_{n}\|_{L^{2}}
\leq \|u_{0}\|_{L^{2}}
\frac{C_{\eps}}{\gamma} e^{\alpha_{*}(\pi+\eps)^{2}/T}
\left( 
\sum_{n=1}^{\infty} e^{-\ldn T} e^{2\eps \sqrt{\ldn}}
\right)^{1/2}
\ .
$$
Since as $T \to 0$ we have
$$\sum_{n=1}^{\infty} e^{-\ldn T} e^{2\eps \sqrt{\ldn}}
\leq e^{2\eps^{2}/T}\sum_{n=1}^{\infty} e^{-\ldn T/2}
\sim e^{2\eps^{2}/T} (T/2)^{-1/2}\Gamma(1/2)/2
\ll C_{\eps}' e^{3\eps^{2}/T} 
\ ,
$$
this implies 
$\displaystyle
\|g\|_{L^{2}(-T/2,T/2)}\leq C_{\alpha}e^{\alpha \pi^{2}/T }\|u_{0}\|_{L^{2}(0,X)} 
$,
with $\alpha=\alpha_{*}(1+2\eps/\pi)^{2} + 3\eps^{2}/\pi^{2}$
and $C_{\alpha}=C_{\eps}C_{\eps}'/\gamma $.
Since $\alpha\to \alpha_{*}$ as $\eps \to 0$,
this completes the proof that theorem~\ref{theo:G}
implies theorem~\ref{theo:1dParabolic}.

\subsection{Complex analysis multipliers}

In this subsection we shall prove theorem~\ref{theo:G}
by the following classical method in complex analysis
(cf. %\cite{LK71} or 
section~14 in~\cite{Red77} for a concise account with references,
and the two volumes~\cite{Koo}
%~\cite{KI98},~\cite{KII92},
for an extensive monograph on multipliers): 
for all $n\in \mathbb{N}^{*}$ and small $\tau>0$, 
we shall form an infinite product $F_{n}$
normalized by $F_{n}(i\ldn)=1$
with zeros at $i\ldk$ for every positive integer $k\neq n$,
and construct a multiplier $M_{n}$ of exponential type $\tau$
with fast decay at infinity on the real axis 
so that $G_{n}=M_{n}F_{n}$ is in $L^{2}$ on the real axis.
At infinity, 
it is well known that the growth of $z \mapsto \ln |F_{n}(z)|$  
can be bounded from above by a power of $|z|$ 
which is inverse to that of $n\mapsto |i\ldn|\sim n^{2}$ 
(cf. theorem~2.9.5 in~\cite{Boa54})
%a theorem of Borel ensures the equality of 
%the exponential order and the exponent of convergence of the zeros 
%of an infinite product of Weierstrass primary factors): 
we prove that our $\ln F_{n}$ 
is essentially bounded by $z \mapsto \pi \sqrt{|z|}+o(\sqrt{\lambda_{n}})$
where the constant $\pi$ is optimal
(cf. remark~\ref{rem:1d}).
%and we give a sharp control on the growth of $c_{n}$ with respect to $n$.
%(the construction of functions $F_{n}$ in lemma~3 in~\cite{SAI00} 
%applies to more general sequences $\{\lambda_{n}\}_{n\in \mathbb{N}^{*}}$
%but yield a greater constant).
Therefore $M_{n}$
has to be essentially bounded by $C_{n}(\tau)\exp(-\pi \sqrt{|x|})$
on the real axis, for some constant $C_{n}(\tau)>0$.
The key point (as in~\cite{Sei84},
theorem~1 in \cite{Sei86} and theorem~2 in~\cite{SAI00}) 
is to construct a multiplier $M_{n}$ such that $C_{n}(\tau)$
has the smallest growth as $\tau$ tends to $0$.
The following two lemmas give the key 
to the construction of $F_{n}$ and $M_{n}$
respectively.
%An extra difficulty that we have to handle is 
%to minimize the growth with respect to $n\in \mathbb{N}^{*}$.

\begin{lem}
\label{lemF}
Let $\{\lambda_{n}\}_{n\in \mathbb{N}^{*}}$ be a sequence of real numbers 
satisfying $(\ref{eqSpec})$.
For all $\eps>0$
there is a $A_{\eps}>0$ such that, 
for all $n\in \mathbb{N}^{*}$, 
the entire function $f_{n}$
defined by 
$\displaystyle
f_{n}(z)=\prod_{k\neq n} \left( 1-\frac{z}{\ldk}\right)
$ 
satisfies
\begin{eqnarray}
&&
  \ln|f_{n}(z)| \leq (\pi+\eps)\sqrt{|z|} + A_{\eps}
\label{eqf1} \\
&&
  \left| \ln|f_{n}(\ldn)| \right|
\leq \eps \sqrt{\ldn} + A_{\eps}
\label{eqf2} 
\end{eqnarray}
\end{lem}

\begin{proof}
For every $n\in \mathbb{N}^{*}$, 
we introduce the counting function of the sequence 
$\{\ldk\}_{k\in \mathbb{N}^{*}\setminus\{n\}}$:
$$
N_{n}(r)=\#\{ k\in \mathbb{N}^{*}\setminus\{n\}  \ | \ \ldk\leq r \}
\ .
$$
From $(\ref{eqSpec})$
%Since$\{\lambda_{n}\}_{n\in \mathbb{N}^{*}}$ is increasing,
we have $N_{0}-1\leq N_{n}\leq N_{0}$ and $\sqrt{\ldn}=n +\nu + o(1)$.
Since $\ldk\leq r <\lkp$ implies 
$\sqrt{\ldk}-k\leq \sqrt{r}-N_{0}(r)\leq \sqrt{\lkp}-(k+1)+1$,
we deduce $|\sqrt{r}-N_{n}(r)-\nu|\leq 2+o(1)$.
The proof uses the assumption $(\ref{eqSpec})$ 
through the estimates of the increments 
$\Lambda_{n}:=\lnp-\ldn$ and $\Delta_{n}:=\slnp-\sldn$
and their increments:
\begin{eqnarray}
\label{eqld}
&&
\ldn=n^{2}+2\nu n+o(n)\ , \quad
\Lambda_{n}=2n+o(n)\ ,  \quad
\Lambda_{n}-\Lambda_{n-1}=o(n), 
\\ 
\label{eqsld}
&& 
\sldn=n + \nu +o(1)\ , \quad
\Delta_{n}=1+o(1)\ ,  \quad
\Delta_{n}-\Delta_{n-1}=o(1)\ , 
\\ 
\label{eqN}
&&
\forall r\in ]0, \lambda_{1}[, N_{n}(r)=0 \ , \quad  
% N_{n}(\lambda_{1}^{-})=0 \ , \quad 
% N_{n}(\lnm^{+})=N_{n}(\lnp^{-})=n-1 \ , \quad  
%|\sqrt{r}-N_{n}(r)|=O(1)
\exists A>0, \ \forall r, \   |\sqrt{r}-N_{n}(r)|\leq A 
\ .
\end{eqnarray}
We shall use repeatedly that for any real sequence 
$\{r_{n}\}_{n\in \mathbb{N}^{*}}$
\begin{equation}
\label{eqr}
r_{n}=o(1) \Rightarrow 
\left| \ln\left(1+\frac{r_{n}}{1+o(1)}\right) \right|
=|r_{n}|(1+o(1))
\end{equation}

To prove $(\ref{eqf1})$,
we estimate the left hand side in terms of $N_{n}$:
\begin{eqnarray*}
\ln|f_{n}(z)| & \leq & 
\sum_{k\neq n} \ln\left( 1+ \frac{|z|}{\ldk} \right)
=\int_{0}^{\infty} \ln\left( 1+ \frac{|z|}{r} \right) dN_{n}(r)
\\
&=&\int_{0}^{\infty} N_{n}(r)\frac{|z|}{|z|+r} \frac{dr}{r}
=\int_{0}^{\infty} \frac{N_{n}(|z|s)}{1+s} \frac{ds}{s}
\end{eqnarray*}
To estimate this last integral 
we use $(\ref{eqN})$ and the integral computations:
\begin{eqnarray*}
\int_{0}^{\infty} \frac{\sqrt{s}}{1+s} \frac{ds}{s}
=\int_{0}^{\infty} \frac{2 dr}{1+r^{2}}=\pi 
\ , \quad
\int_{\frac{\lambda_{1}}{|z|}}^{\infty} 
\frac{ds}{s(1+s)}
=\left[ \ln\left|\frac{s}{1+s}\right|\right]_{\frac{\lambda_{1}}{|z|}}^{\infty} 
=\ln(1+\frac{|z|}{\lambda_{1}})
\end{eqnarray*}
Thus we obtain 
$\ln|f_{n}(z)|\leq 
\pi\sqrt{|z|} + A\ln(1+\frac{|z|}{\lambda_{1}})$,
so that, for all $\eps>0$ there is a $A_{\eps}'>0$ such that   
$ \ln|f_{n}(z)|\leq (\pi+\eps)\sqrt{|z|} + A_{\eps}'$.

To prove $(\ref{eqf2})$,
we estimate the left hand side in terms of $N_{n}$:
\begin{eqnarray*}
\ln|f_{n}(\ldn)| 
& = &
\sum_{k< n} \ln\left( \frac{\ldn}{\ldk} - 1\right)
+ \sum_{k> n} \ln\left( 1- \frac{\ldn}{\ldk} \right) \\
{}& = &
\int_{\lambda_{1}^{-}}^{\lnm^{+}} \ln\left( \frac{\ldn}{r} - 1\right)dN_{n}(r)
+ \int_{\lnp^{-}}^{\infty} \ln\left( 1-\frac{\ldn}{r} \right)dN_{n}(r)
\end{eqnarray*}
Integrating by parts yields $\ln|f_{n}(\ldn)|=I_{n}+B_{n}$ with 
\begin{eqnarray*}
I_{n}&=&\int_{\lambda_{1}^{-}}^{\lnm^{+}} 
N_{n}(r) \frac{\ldn}{\ldn-r} \frac{dr}{r}
+ \int_{\lnp^{-}}^{\infty}
N_{n}(r) \frac{\ldn}{\ldn-r} \frac{dr}{r}
\\
B_{n}&=&
\left[N_{n}(r)\ln\left( \frac{\ldn}{r} - 1\right)
\right]_{\lambda_{1}^{-}}^{\lnm^{+}}
+ 
\left[N_{n}(r)\ln\left( 1-\frac{\ldn}{r} \right)
\right]_{\lnp^{-}}^{\infty}
\end{eqnarray*}
To estimate the boundary term $B_{n}$,
we first simplify its expression using  
$N_{n}(\lambda_{1}^{-})=0$
and $N_{n}(\lnm^{+})=N_{n}(\lnp^{-})=n-1$, 
then we sort out the increments $\Lambda_{n}=\lnp-\ldn$,
and finally we use $(\ref{eqld})$ and $(\ref{eqr})$:
\begin{eqnarray*}
B_{n}&=&
(n-1)\left[ 
\ln\left( \frac{\ldn}{\lnm} - 1\right)
-
\ln\left( 1- \frac{\ldn}{\lnp} \right)
\right]
\\
&=&(n-1)\left[ 
\ln\left( 1-\frac{\Lambda_{n}-\Lambda_{n-1}}{\Lambda_{n}} \right)
+\ln\left( 1+\frac{\Lambda_{n}+\Lambda_{n-1}}{\lnm} \right)
\right]
\\
&=&(n-1)\left[
\frac{o(n)}{2n}(1+o(1))
+\frac{4n+o(n)}{n^{2}}(1+o(1))
\right]
=o(1)
\ .
\end{eqnarray*}
Now we estimate the integral term $I_{n}$.
%The change of variable $r=\ldn s$ yields:
% \begin{eqnarray*}
% I_{n}&=&\int_{\frac{\lambda_{1}^{-}}{\ldn}}^{\frac{\lnm^{+}}{\ldn}} 
% N_{n}(\ldn s) \frac{ds}{s(1-s)} 
% + \int_{\frac{\lnp^{-}}{\ldn}}^{\infty}
% N_{n}(\ldn s) \frac{ds}{s(1-s)} 
% \end{eqnarray*}
Performing the change of variable $r=\ldn s$ 
and using $(\ref{eqN})$ yields:
$|I_{n}-\sqrt{\ldn}J_{n}|\leq A K_{n}$
with 
\begin{eqnarray*}
J_{n}&=&
\int_{\frac{\lambda_{1}^{-}}{\ldn}}^{\frac{\lnm^{+}}{\ldn}} 
\frac{ds}{(1-s)\sqrt{s}}
+ \int_{\frac{\lnp^{-}}{\ldn}}^{\infty}
\frac{ds}{(1-s)\sqrt{s}}
\\
K_{n}&=&\int_{\frac{\lambda_{1}^{-}}{\ldn}}^{\frac{\lnm^{+}}{\ldn}} 
\frac{ds}{(1-s)s}
+ \int_{\frac{\lnp^{-}}{\ldn}}^{\infty}
\frac{ds}{(s-1)s}
\end{eqnarray*}
The term $K_{n}$ is readily computed and estimated using $(\ref{eqld})$:
\begin{eqnarray*}
K_{n}&=&
\left[\ln\frac{s}{1-s}
\right]_{\frac{\lambda_{1}^{-}}{\ldn}}^{\frac{\lnm^{+}}{\ldn}} 
+ 
\left[\ln\frac{s-1}{s}
\right]_{\frac{\lnp^{-}}{\ldn}}^{\infty}
=\ln\frac{\lnp}{\Lambda_{n}}
+\ln\frac{\lnm}{\Lambda_{n-1}}
+\ln\left(\ldn(\frac{1}{\lambda_{1}}+\frac{1}{\ldn})\right)
%(\lambda_{1}^{-1}-\ldn^{-1})\right)
%+2\ln\left(\sldn\sqrt{\lambda_{1}^{-1}-\ldn^{-1}}\right)
\\
&=& 
2\ln\frac{n^{2}+O(n)}{2n+o(n)}
+2\ln\sldn + O(1)
=o(\sldn)
\ .
\end{eqnarray*}
We compute $J_{n}$ after a change of variable, 
and estimate it by $(\ref{eqsld})$ and $(\ref{eqr})$
after sorting out the increments $\Delta_{n}=\slnp-\sldn$:
\begin{eqnarray*}
J_{n}
&=&
\int_{\frac{\sldo^{-}}{\sldn}}^{\frac{\slnm^{+}}{\sldn}} 
\frac{2dr}{r^{2}-1}
+ \int_{\frac{\slnp^{-}}{\sldn}}^{\infty}
\frac{2dr}{r^{2}-1}
=
\left[
\ln{\frac{1-r}{r+1}}
\right]_{\frac{\sldo^{-}}{\sldn}}^{\frac{\slnm^{+}}{\sldn}} 
+\left[
\ln{\frac{r-1}{r+1}}
\right]_{\frac{\slnp^{-}}{\sldn}}^{\infty}
\\
&=&
\ln\left(\frac{\Delta_{n-1}}{\Delta_{n}}\right)
+\ln\left(\frac{\slnp+\sldn}{\sldn+\slnm}\right)
-\ln\left(\frac{\sldn-\sldo}{\sldn+\sldo}\right)
\\
&=&
\ln\left(1-\frac{\Delta_{n}-\Delta_{n-1}}{\Delta_{n}}\right)
+\ln\left(1+\frac{\Delta_{n}+\Delta_{n-1}}{\sldn+\slnm}\right)
-\ln\left(1-\frac{2\sldo}{\sldn+\sldo}\right)
\\
% &=&
% \ln\left(1+\frac{o(1)}{1+o(1)}\right) 
% +\ln\left(1+\frac{2+o(1)}{2n(1+o(1))}\right) 
% -\ln\left(1-\frac{O(1)}{n(1+o(1))}\right)
% \\
&=&
o(1)(1+o(1))
+\frac{2+o(1)}{2n}(1+o(1))
+\frac{O(1)}{n}(1+o(1))
%= o(1) + O(n^{-1})+ O(n^{-1})
=o(1)
\ .
% J_{n} &=&
% \ln\left(1+\frac{\Delta_{n}+\Delta_{n-1}}{\Delta_{n-1}}\right)
% -\ln\left(1-\frac{2\sldo}{\sldn+\sldo}\right)
% +\ln\left(1-\frac{\Delta_{n}-\Delta_{n-1}}{\Delta_{n}}\right)
% J_{n} &=&
%\ln\left(\frac{\Delta_{n}}{\slnp+\sldn}\right)
%-\ln\left(\frac{1-\sldo/\sldn}{1+\sldo/\sldn}\right)
%+\ln\left(1-\frac{\Delta_{n}-\Delta_{n-1}}{\Delta_{n}}\right)
\end{eqnarray*}
Plugging the estimates $K_{n}=o(\sldn)$ and $J_{n}=o(1)$ 
into $|I_{n}-\sqrt{\ldn}J_{n}|\leq A K_{n}$
yields $I_{n}=o(\sldn)$.
Plugging this estimate and $B_{n}=o(1)$
into $\ln|f_{n}(\ldn)|=I_{n}+B_{n}$
yields $\ln|f_{n}(\ldn)|=o(\sldn)$,
which completes the proof of $(\ref{eqf2})$.
\end{proof}

\begin{lem}
\label{lemM}
Let $\alpha_{*}$ be defined by $(\ref{eqalphastar})$.
For all $d>0$ there is a $D>0$
such that for all $\tau>0$, %$\tau\in ]0,1]$, 
there is an even entire function $M$ of exponential type 
%(lower or equal to) 
$\tau$ satisfying: 
%such that:
$M(0)=1$ and
%, $M$ is even 
%and satisfies
\begin{equation}
\label{eqM}
%M(0)=1\ \mbox{, } \quad 
\forall x>0, \quad \ln|M(x)|\leq \frac{\alpha_{*}d^{2}}{2\tau}+D-d\sqrt{x} 
\quad \mbox{and} \quad |M(ix)|\geq 1
\ .
\end{equation}
\end{lem}
\begin{proof}
Following Ingham and many others since 1934
(cf. section~14 in~\cite{Red77} for theorems and references)
we seek a multiplier $M$ of small exponential type 
decaying rapidly along the real axis 
in the following form: 
\begin{equation}
\label{eqMIngham}
M(z)=\prod_{n\in \mathbb{N}}\sinc\left(\frac{ z}{a_{n}}\right)
\quad \mbox{where} \
\sinc(0)=1,\ \forall z\in \mathbb{C}^{*},\ \sinc(z)=\frac{\sin(z)}{z}
\end{equation}
and where $\{a_{n}\}_{n\in \mathbb{N}}$ is a non decreasing %increasing 
sequence of positive real numbers such that 
$\tau_{M}=\sum_{n\in \mathbb{N}}\frac{1}{a_{n}}<\infty$.
Since the cardinal sine function $\sinc$ is an even entire function
of exponential type $1$ 
satisfying $\sinc(0)=1$ and $\sinc(ix)=\sinh(x)/x\geq 1$ for all $x>0$, 
$(\ref{eqMIngham})$ defines an even entire function $M$ 
of exponential type $\tau_{M}$ 
satisfying $M(0)=1$ and $|M(ix)|\geq 1$ for all $x>0$.

We define $\{a_{n}\}_{n\in \mathbb{N}}$
by the slope $A$ of its counting function $N$
and its first term $a_{0}$ (to be chosen large enough):
$$
N(r):=\sum_{|a_{n}|\leq r}1 
=[A\sqrt{u}] 
\ \mbox{for} \ r\geq 2
\quad \mbox{and} \quad 
a_{0}\geq A^{-2}
\ ,
$$
where $[x]$ denotes as usual 
the greatest integer smaller or equal to the real number $x$.
The exponential type $\tau_{M}$ of $M$ is easily
bounded from above by $\tau=2A/\sqrt{a_{0}}$:
$$
\tau_{M}:=\sum_{n\in \mathbb{N}}\frac{1}{a_{n}}
=\int_{0}^{\infty}\frac{dN(r)}{r}
=\int_{0}^{\infty}\frac{N(r)}{r^{2}}dr
\leq \int_{a_{0}}^{\infty}\frac{A\sqrt{r}}{r^{2}}dr
=\frac{2A}{\sqrt{a_{0}}}
=: \tau
\ ,
$$
and we are left with estimating the decay of:
\begin{equation}
\label{eqlnM}
\ln|M(x)|
=\int_{a_{0}^{-}}^{\infty}f\left(\frac{x}{r}\right)\, dN(r)
\quad \mbox{where} \
f(\theta)=\ln \sinc(\theta)=\ln \frac{\sin(\theta)}{\theta}
\ .
\end{equation}
We shall choose $A$ such that, for all  $a_{0}\geq A^{-2}$,
$\ln|M(x)|\leq -d\sqrt{x}+ O(1)$ as $x\to +\infty$,
and then prove that:
$\ln|M(x)|\leq \alpha_{*}d^{2}/(2\tau)-d\sqrt{x}+ O(1)$ as $\tau\to 0$
(equivalently $a_{0}\to +\infty$)
uniformly in $x>0$.

\medskip

For $x>a_{0}$ we take advantage of the boundedness of sine 
through the estimate $f(\theta)\leq -\ln|\theta|$ for $|\theta|\leq 1$,
by splitting the integral in $(\ref{eqlnM})$ into the two terms:
\begin{eqnarray*}
I&=& 
\int_{a_{0}^{-}}^{x}f\left(\frac{x}{r}\right)\, dN(r)
\leq \int_{a_{0}^{-}}^{x}\ln \left|\frac{r}{x}\right|\, dN(r)
= -\int_{a_{0}}^{x}N(r)\frac{dr}{r}
\\
J&=&
\int_{x}^{\infty}f\left(\frac{x}{r}\right)\, dN(r)
%=\int_{1}^{\infty}f'\left(\frac{1}{s}\right)
%\frac{N(xs)}{s}\frac{ds}{s}-f(1)N(x)
=\int_{0}^{1}f'(\theta)
N\left(\frac{x}{\theta}\right)
\, d\theta
-f(1)N(x)
\end{eqnarray*}
where right hand sides were integrated by parts
and $\theta=x/r$.
Now we plug in the basic estimate on $N$: 
$A\sqrt{r}-1\leq N(r)\leq A\sqrt{r}$ for $r\geq 2$.
The first term is now estimated by 
\begin{equation}
\label{eqlnM1}
I\leq -A\int_{a_{0}}^{x}\frac{dr}{\sqrt{r}}+
\int_{a_{0}}^{x}\frac{dr}{r}
= -2A\left(\sqrt{x}-\sqrt{a_{0}}\right)
+\ln x - \ln a_{0}
\ .
\end{equation}
To estimate the second term, we first observe 
that the Hadamard factorization of the cardinal sine function
$\sinc(\pi z)=\prod_{n\in \mathbb{N}^{*}}\left(1-\frac{z^{2}}{n^{2}}\right)$
and the Taylor expansion of the logarithm at $1$ imply:
$$
f(\theta)=-\sum_{k\in \mathbb{N}^{*}}
\frac{\zeta(2k)}{k}\left(\frac{\theta}{\pi}\right)^{2k}
\ \mbox{for} \ |\theta|< 1
\quad \mbox{, where} \
\zeta(s)=\sum_{n\in \mathbb{N}^{*}}\frac{1}{n^{s}}
\ .
$$
The second term is now estimated by 
\begin{eqnarray}
\label{eqlnM2}
J &\leq &
\int_{0}^{1}f'(\theta)
\left(\frac{A\sqrt{x}}{\sqrt{\theta}}-1\right)d\theta
-f(1)A\sqrt{x}
\\
\nonumber
&=&
A\sqrt{x}\left(
\int_{0}^{1}f'(\theta)\frac{d\theta}{\sqrt{\theta}}  -f(1)
\right)
-f(1)
% =
% A\sqrt{x}\int_{0}^{1}f'(\theta)
% \frac{d\theta}{\sqrt{\theta}}
% -f(1)(A\sqrt{x}+1)
\\
\nonumber
&=& 
-A\sqrt{x}
\sum_{k\in \mathbb{N}^{*}} 
\left(
\frac{2k}{2k-\frac{1}{2}}-1
\right)
\frac{\zeta(2k)}{k\pi^{2k}}
-f(1)
=-A\Sigma^{*}\sqrt{x}-f(1)
\ ,
\end{eqnarray}
where the series for $f$ was differentiated, multiplied and integrated 
term by term, 
and $\Sigma^{*}=\sum_{k\in \mathbb{N}^{*}} 
\frac{1}{k(4k-1)}
\frac{\zeta(2k)}{\pi^{2k}}$.
Putting $(\ref{eqlnM1})$ and $(\ref{eqlnM2})$ together yields:
\begin{equation*}
%\label{eqlnM3}
\forall x>a_{0}\ , \
  \ln|M(x)|\leq -(2+\Sigma^{*})A\sqrt{x}
+\ln x -f(1)
+ 2A\sqrt{a_{0}}
\ .
\end{equation*}
so that, for all $d>(2+\Sigma^{*})A$ there is a $D_{1}$ such that 
\begin{equation}
\label{eqlnM3}
\forall d>(2+\Sigma^{*})A, \exists D_{1}>0, 
\forall x>a_{0}, \
  \ln|M(x)|\leq 2A\sqrt{a_{0}}-d\sqrt{x}+D_{1} 
\ .
\end{equation}

\medskip

Since $|\sinc|$ is bounded by $1$:
for all $x$, $\ln|M(x)|\leq 0$. 
Moreover  $d>2A$, 
so that $(\ref{eqlnM3})$ implies 
\begin{equation}
\label{eqlnM4}
\forall a_{0}\geq A^{-2}, \forall x>0, \
  \ln|M(x)|\leq d\sqrt{a_{0}}-d\sqrt{x}+D_{1} 
\ .
\end{equation}
Since $d>(2+\Sigma^{*})A$ and 
$\tau=2A/\sqrt{a_{0}}$,
this proves: % for all $\tau \leq 2A^{2}$:
\begin{equation}
\label{eqMM1}
\forall \tau\leq 2A^{2},\forall x>0, \quad 
\ln|M(x)|\leq \frac{\alpha_{1}d^{2}}{2\tau}-d\sqrt{x}+D_{1} 
\end{equation}
with $\alpha_{1}=4/(2+\Sigma_{*})$.

\medskip

For $x<a_{0}$, we can also use the better estimate:
\begin{eqnarray}
\label{eqlnM5}
\ln|M(x)| &\leq &
\int_{a_{0}}^{\infty}f\left(\frac{x}{r}\right)\, dN(r)
=\int_{0}^{x/a_{0}}f'(\theta)
N\left(\frac{x}{\theta}\right)
\, d\theta
\\
\nonumber
&\leq & 
A\sqrt{x}
\int_{0}^{x/a_{0}}f'(\theta)\frac{d\theta}{\sqrt{\theta}}  
-f\left(\frac{x}{a_{0}}\right)
\\
\nonumber
&\leq &-A\sqrt{a_{0}}
\sum_{k\in \mathbb{N}^{*}} 
\frac{4k\zeta(2k)}{k(4k-1)}
\left(\frac{x}{a_{0}\pi}\right)^{2k}
-f\left(1\right)
\end{eqnarray}
If we keep only the first term (i.e. $k=1$) of the series 
in $(\ref{eqlnM3})$ and $(\ref{eqlnM5})$, we get that
for all $d>(2+\frac{1}{3}\frac{\zeta(2)}{\pi^{2}})A$ 
there is a $D_{2}$ such that :
\begin{eqnarray}
\label{eqlnM6}
\forall x>a_{0}, \
  \ln|M(x)|
&\leq& 2A\sqrt{a_{0}}-d\sqrt{x}+D_{2} 
\\
\nonumber
\forall x<a_{0}, \
\ln|M(x)| &\leq &
-A\sqrt{a_{0}}
\frac{4\zeta(2)}{3\pi^{2}}
\left(\frac{x}{a_{0}}\right)^{2}
-f\left(1\right)
\ .
\end{eqnarray}
Now, for all $x<a_{0}$
\begin{eqnarray*}
%  \forall x<a_{0}, \
\ln|M(x)|- 2A\sqrt{a_{0}}+d\sqrt{x}
&\leq &
A\sqrt{a_{0}}F\left(\frac{x}{a_{0}}\right)
% \left[
% -2 + (2+\eps)\sqrt{X}+\frac{1}{3}\frac{\zeta(2)}{\pi^{2}}(\sqrt{X}-4X^{2})
% \right]
\end{eqnarray*}
with 
$F(X)=-2 + (2+\eps)\sqrt{X}+\frac{1}{3}\frac{\zeta(2)}{\pi^{2}}(\sqrt{X}-4X^{2})
=-2+(37/18+\eps)\sqrt{X}-2X^{2}/9$
and $\eps=d/A-(2+\frac{1}{3}\frac{\zeta(2)}{\pi^{2}})>0$.
Since $F$ is increasing on $[0,1]$ 
and $F(1)=\eps-1/6$,
choosing $A$ so that $\eps<1/6$, 
yields that 
$
\ln|M(x)|- 2A\sqrt{a_{0}}+d\sqrt{x}\leq 0
$,
for all $x<a_{0}$.
Together with $(\ref{eqlnM6})$, this proves 
\begin{eqnarray}
%\label{eqlnM6}
\forall x>0, \
  \ln|M(x)|
&\leq& 2A\sqrt{a_{0}}-d\sqrt{x}+D_{2} 
\end{eqnarray}
Since 
$d>(2+\frac{1}{3}\frac{\zeta(2)}{\pi^{2}})A=37A/18$ and 
$\tau=2A/\sqrt{a_{0}}$,
this proves: % for all $\tau \leq 2A^{2}$:
\begin{equation}
\label{eqMM2}
\forall \tau\leq 2A^{2},\forall x>0, \quad 
\ln|M(x)|\leq \frac{\alpha_{2}d^{2}}{2\tau}-d\sqrt{x}+D_{2} 
\end{equation}
with $\alpha_{2}=2(36/37)^{2}$.

\medskip 

Equations $(\ref{eqMM1})$ and $(\ref{eqMM2})$ 
complete the proof of the lemma~\ref{lemM}
with $\alpha_{*}=\min\{\alpha_{1},\alpha_{2}\}$.
Since we have checked on a computer that $\alpha_{1}>\alpha_{2}$,
we decided to state the lemma with $\alpha_{*}=\alpha_{2}$,
i.e. $(\ref{eqalphastar})$. 
% As in~\cite{Sei84} and~\cite{Sei86}, 
% we shall use the two estimates
% on the sine cardinal function $\sinc$:
% \begin{equation}
% \forall \theta \in \mathbb{R}, \ %\quad %[-1,1], \
% 0<\sinc(\theta)\leq \exp\left(-\frac{\theta^{2}}{6}\right) 
% %e^{-\frac{\theta^{2}}{6}}
% \ \mbox{for} \ |\theta|\leq 1
% \ , \ %\quad \mbox{and} \quad 
% |\sinc(\theta)|\leq \frac{1}{|\theta|}
% \ \mbox{for} \ |\theta|\geq 1
% %\ .
% \end{equation}
% (the first one is obtained by comparing the Taylor series expansions).
\end{proof}

To prove theorem~\ref{theo:G},
we use lemmas~\ref{lemF} and~\ref{lemM}
with $d=\pi+2\eps$  
and define :
%$G_{n}$ as:
\begin{equation*}
G_{n}=F_{n}M_{n}
\quad \mbox{with} \quad 
F_{n}(z)=f_{n}(-iz)/f_{n}(\ldn)
\quad \mbox{and} \quad
M_{n}(z)=M(z)/M(i\ldn)
\ .
\end{equation*}

Thanks to lemma~\ref{lemF},
the entire function $F_{n}$
satisfies
\begin{eqnarray}
&&
F_{n}(i\ldn)=1 
\quad \mbox{ and  } \quad
\forall k\in \mathbb{N}^{*} ,\, k\neq n ,\,  
F_{n}(i\ldk)=0 , 
\label{eqF1} \\
&&
  \ln|F_{n}(z)| \leq (\pi+\eps)\sqrt{|z|} + \eps \sqrt{\ldn} + 2A_{\eps}
\label{eqF2} 
\end{eqnarray}
where $(\ref{eqF1})$ is an obvious consequence of 
the definitions of $f_{n}$ and $F_{n}$, 
and $(\ref{eqF2})$ is a consequence of 
the estimates $(\ref{eqf1})$ and $(\ref{eqf2})$.

Thanks to lemma~\ref{lemM}, %with $d=d_{\eps}=\pi+2\eps$,
there is a $D_{\eps}>0$
such that the entire function $M_{n}$ is of exponential type $\tau$
and satisfies 
\begin{eqnarray}
&&
M_{n}(i\ldn)=1 
\label{eqM1} \\
&&
\forall x\in \mathbb{R}, \quad 
\ln|M_{n}(x)|\leq \frac{\alpha_{*}d^{2}}{2\tau}+D_{\eps}-d\sqrt{|x|}
\label{eqM2} 
\end{eqnarray}
where $(\ref{eqM1})$ is an obvious consequence of the definitions 
of $M$ and $M_{n}$, 
and $(\ref{eqM2})$ is a consequence of $(\ref{eqM})$ since $M$ is even. 

The entire function $G_{n}$ 
has the same exponential type as $M_{n}$
since $(\ref{eqF2})$ 
implies that the exponential type of $F_{n}$ is $0$.
Hence $(\ref{eqG1})$ holds.
Putting $(\ref{eqF1})$ and $(\ref{eqM1})$ together yields $(\ref{eqG2})$.
Since $d=\pi+2\eps$, 
$(\ref{eqF2})$ and $(\ref{eqM2})$ imply
$$
 \forall x\in \mathbb{R}, \quad 
\ln|G_{n}(x)| \leq 
D_{\eps}+2A_{\eps}-\eps\sqrt{|x|} 
+\eps \sqrt{\ldn} + \frac{\alpha_{*}d^{2}}{2\tau} 
\ .
$$
Hence $(\ref{eqG3})$ holds
with 
$\displaystyle
C_{\eps}=e^{D_{\eps}+2A_{\eps}}\left(
\int_{-\infty}^{+\infty} e^{ -2\eps\sqrt{|x|} } \, dx
\right)^{1/2} $.
Theorem~\ref{theo:G} is proved.
%This completes the proof of theorem~\ref{theo:G}.

\begin{rem}
\label{rem:1d}
Under assumption $(\ref{eqSpec})$,
lemma~3 in \cite{SAI00} (which applies to much more general sequences)
proves that %the function 
$\displaystyle
F_{n}(z)=\prod_{k\neq n} \left[ 1- 
\left( \frac{z-\ldn}{\ldk-\ldn} \right)^{2}
\right]
$ 
satisfies $(\ref{eqF1})$ and 
$\ln|F_{n}(\ldn + z)| \leq 2\pi\sqrt{|z|}$,
hence
$\ln|F_{n}(z)| \leq 2\pi\sqrt{|z|}+O(\sqrt{\ldn})$.
In $(\ref{eqF2})$,
the estimate $O(\sqrt{\ldn})$ improves to $o(\sqrt{\ldn})$ 
and the constant $2\pi$ improves to the optimal $\pi$
(optimality can be deduced from theorem~4.1.1 in~\cite{Boa54}).

Seidman obtained lemma~\ref{lemM} for
$\alpha_{*}=\beta_{*}$
with $\beta_{*}\approx 42.86$
in the proof of Theorem~3.1 in \cite{Sei84}.
His later Theorem~1 in \cite{Sei86}
improves the rate to $\alpha_{*}=2\beta_{*}$ 
with $\beta_{*}\approx 4.17$.
Theorem~2 in \cite{SAI00}, 
which applies to much more general spectral sequences,
yields lemma~\ref{lemM} for $\alpha_{*}=24$.
The argument used in section~\ref{sec:lb}
can be used to prove that lemma~\ref{lemM}
does not hold for $\alpha_{*}<1/4$.
It would be interesting to determine the smallest value of $\alpha_{*}$
for which it holds.
% We do not claim that the value $\alpha_{*}$
% defined by $(\ref{eqalphastar})$ 
% in lemma~\ref{lemM} is optimal.
% It falls short of 
% %Yet we have not been able to reach 
% the value $\alpha_{*}=1/4$
% which is a lower bound 
% that would prove a one dimensionnal version of the 
% conjecture stated at the end of the first subsection 
% in section~\ref{sec:res}.
\end{rem}

%%%%%%%%%%%%%%%%%%%%%%%%%%%%%%%%%%%%%%%%%%%%%%%%%%%%%%%%%%%%
%% TRANSMUTATION
%%%%%%%%%%%%%%%%%%%%%%%%%%%%%%%%%%%%%%%%%%%%%%%%%%%%%%%%%%%%

\section{Upper bound under the geodesics condition}
\label{sec:transmut}

In this section we prove 
theorem~\ref{theo:ub}
in three steps.
$ \mathcal{D}'(\mathcal{O})$ denotes the space of 
distributions on the open set $\mathcal{O}$
endowed with the weak topology 
and $\mathcal{M}(\mathcal{O})$ denotes the subspace of 
Radon measures on $\mathcal{O}$.

\subsection{The segment controlled at both ends}

In a first step we prove that the upper bound 
for the null-controllability cost 
of the heat equation on the segment $[0,L]$ 
controlled at one end
is the same as the null-controllability cost 
of the heat equation on the twofold segment $[-L,L]$ 
controlled at both ends.

Given a time $T>0$ and a length $L>0$,
we denote by $D$ (respectively $N$) 
some continuous operator from $L^{2}(0,L)$  to $L^{2}(0,T)$
allowing to control to zero in time $T$ 
the heat equation on $[0,L]$
with zero Dirichlet (respectively Neumann) condition at $0$
by a Dirichlet control at $L$.
More precisely, for all $u_{0}\in L^{2}(0,L)$
the solution $u\in C^{0}([0,\infty),L^{2}(0,L))$,
denoted by $u=S_{D}u_{0}$ (respectively $u=S_{N}u_{0}$), 
of the Cauchy problem in theorem~\ref{eqHeat1d} 
with $B=1$ (respectively $B=\partial_{s}$)
and $g=Du_{0}$ (respectively $g=Nu_{0}$)
satisfies $u=0$ at $t=T$.
\begin{prop}
\label{prop:twofold}
For any time $T>0$ and any length $L>0$,
there is a continuous operator $K$ 
from $L^{2}(-L,L)$  to $L^{2}(0,T)^{2}$
allowing to control to zero in time $T$ 
the heat equation on $[-L,L]$
by Dirichlet controls at both ends
at the same cost as $D$ and $N$,
i.e. for all $v_{0}\in L^{2}(-L,L)$
the solution \linebreak $v\in C^{0}([0,\infty),L^{2}(-L,L))$ of~:
%denoted by $v=Sv_{0}$, of~: 
\begin{equation} 
\label{eqHeattwofold}
\partial_{t}v - \partial_{s}^{2} v=0
\quad {\rm in}\ ]0,T[\times ]-L,L[ ,\quad 
(v\res{s=-L} ,v\res{s=L})=Kv_{0} ,\quad
v\res{t=0} = v_{0} 
\end{equation}
satisfies $v=0$ at $t=T$ and 
$\displaystyle
\|K\|\leq\sup(\|D\|,\|N\|)$.
\end{prop}
\begin{proof}
Given $v_{0}\in L^{2}(-L,L)$,
we decompose it in odd and even parts~: $v_{0}=v_{0,odd}+v_{0,even}$.
We denote by $u_{0,odd}$ and $u_{0,even}$
the restrictions of $v_{0,odd}$ and $v_{0,even}$ to $[0,L]$,
We denote by $f=Du_{0,odd}$ and $g=Nu_{0,even}$ 
the corresponding controls.
We denote by $u_{odd}=S_{D}u_{0,odd}$ and $u_{even}=S_{N}u_{0,even}$ 
the corresponding solutions.

We define $v\in L^{2}([0,T]\times[-L,L])$ by 
$v(t,\pm s)=u_{even}(t,s) \pm u_{odd}(t,s)$ for $s\geq 0$.
Since 
$$
(\partial_{t}-\partial_{s}^{2})u_{even}=
(\partial_{t}-\partial_{s}^{2})u_{odd}= 0
\mbox{ in }
\mathcal{D}'(]0,T[\times]0,L[)
\ ,
$$
we have, denoting the Dirac mass at $s=0$ by $\delta_{s}\in \mathcal{D}'(\mathbb{R})$,
$$(\partial_{t}-\partial_{s}^{2})v
=
2u_{odd}(t,0) \otimes \delta_{s}'(0)
+
2\partial_{s}u_{even}(t,0) \otimes \delta_{s}(0)
\ .
$$
But $u_{odd}(t,0)=\partial_{s}u_{even}(t,0)=0$
by the definition of $D$ and $N$.
Hence  $(\partial_{t}-\partial_{s}^{2})v=0$.
Moreover $v(0,s)=v_{0}(s)$, 
$v(T,s)=0$, 
$v(t,L)=g(t)+f(t)$, 
$v(t,-L)=g(t)-f(t)$.
Therefore, setting $Kv_{0}=(g-f,g+f)$ yields 
an operator $K$ satisfying 
the null-controllability property required.

To finish the proof we estimate its cost $\|K\|$.
Taking the Euclidean norm for $Kv_{0}=(g-f,g+f)$,
we have $\|Kv_{0}\|_{L^{2}(0,T)^{2}}^{2}
= 2\|f\|_{L^{2}(0,T)}^{2}+ 2\|g\|_{L^{2}(0,T)}^{2}$.
Since $f=Du_{0,odd}$ and $g=Nu_{0,even}$,
setting $C=\sup(\|D\|,\|N\|)$ we have 
\begin{eqnarray}
\|Kv_{0}\|_{L^{2}(0,T)^{2}}^{2}
&\leq&
2C^{2}\left( 
\|u_{0,odd}\|_{L^{2}(0,L)}^{2}
+\|u_{0,even}\|_{L^{2}(0,L)}^{2}
\right)
\label{eqKv0} 
\end{eqnarray}
Moreover, 
since $2u_{0,odd}(s)=v_{0}(s)-v_{0}(-s)$ 
and
$2u_{0,even}(s)=v_{0}(s)+v_{0}(-s)$ 
for $s\in [0,L]$,
we have 
\begin{eqnarray}
  \|2u_{0,odd}\|_{L^{2}(0,L)}^{2} &=& 
\|v_{0}\|_{L^{2}(-L,L)}^{2} 
-2 \int_{0}^{L} v_{0}(s)v_{0}(-s) \, ds \label{eqodd} \\
  \|2u_{0,even}\|_{L^{2}(0,L)}^{2} &=& 
\|v_{0}\|_{L^{2}(-L,L)}^{2} 
+2 \int_{0}^{L} v_{0}(s)v_{0}(-s)\, ds \label{eqeven}
\ .
\end{eqnarray}
Equations (\ref{eqKv0}), (\ref{eqodd}) and (\ref{eqeven})
imply 
$\dis
\|Kv_{0}\|_{L^{2}(0,T)^{2}}\leq C \|v_{0}\|_{L^{2}(-L,L)} 
$. 
\end{proof}

\subsection{The fundamental controlled solution}

In a second step we %use theorem~\ref{theo:1d} to 
construct a ``fundamental controlled solution'' $v$
of the heat equation on the segment controlled 
by Dirichlet conditions at both ends.
%use theorem~\ref{theo:1d} to 

\begin{prop}
\label{prop:fundcontrol}
If theorem~\ref{theo:1d} holds for some rate $\alpha_{*}$,
then for any $\alpha > \alpha_{*}$,
there exists $A>0$ such that 
for all $L>0$ and $T\in\, ]0,\inf(\pi/2,L)^{2}]$
%, $L\geq \sqrt{T}$  
%and $u_{0}\in L^{2}(0,L)$
there is a $v\in C^{0}([0,T], \mathcal{M}(]-L,L[))$
satisfying 
\begin{eqnarray}
\partial_{t}v - \partial_{s}^{2}v  =  0
\quad {\rm in }\ \mathcal{D}'(]0,T[\times ]-L,L[)\ , 
\label{eqv1} \\
%v\in C^{0}(]0,T], L^{2}(]-L,L[)), \ 
v\res{t=0}  =  \delta \quad {\rm and }\quad v\res{t=T}  =  0 \ ,
\label{eqv2} \\
\|v\|_{L^{2}(]0,T[\times ]-L,L[)}
\leq Ae^{\alpha L^{2}/T } \ .
\label{eqv3}
\end{eqnarray}
\end{prop}    
We shall sometimes refer to a function $v$ 
satisfying the above requirements 
as a fundamental controlled solution
on $]0,T[\times ]-L,L[$ at cost $(A,\alpha)$.
\begin{proof}
We first reduce the problem to the case $L=\pi/2$
using the rescaling 
$(t,s)\mapsto (\sigma^{2} t, \sigma s)$, $\sigma>0$
with $\sigma=\pi/(2L)$.
Given $L>0$ and $T\in\, ]0,\inf(\pi/2,L)^{2}]$,
%$T\in ]0,1]$ and $L\geq \sqrt{T}$, 
we set $\Lt=\pi/2$ and $\Tt=\sigma^{2}T\in ]0,\Lt^{2}]$. 
Let $\vt$ be a fundamental controlled solution
on $]0,\Tt[\times ]-\Lt,\Lt[$ 
at cost $(\tilde{A},\tilde{\alpha})$.
Setting $v(t,s)=\sigma \vt(\sigma^{2}t,\sigma s)$
defines a fundamental controlled solution $v$
on $]0,T[\times ]-L,L[$ 
at cost $(\tilde{A}/\sqrt{\sigma},\tilde{\alpha})$.
Since $T\leq\Lt^{2}$, we have 
$\tilde{A}/\sqrt{\sigma} 
\leq \tilde{A} (L^{2}/T)^{1/4}$.
Hence for all $\alpha>\tilde{\alpha}$
there is an $A>0$ such that 
$v$ is also a fundamental controlled solution 
on $]0,T[\times ]-L,L[$ at cost $(A,\alpha)$. 
Therefore, it is enough to prove proposition~\ref{prop:fundcontrol}
in the particular case $L=\pi/2$.

We assume theorem~\ref{theo:1d} holds for some rate $\alpha_{*}$.
Let $\at > \at_{*} > \alpha_{*}$, 
$L=\Lt=\pi/2$  and $\Tt\in ]0,\Lt^{2}]$  
be fixed from now on.
We set $\alpha=(1-\eps)\at_{*}$ 
and $T=(1-\eps)\Tt$
where $\eps\in ]0,1[$ is chosen 
such that $\alpha>\alpha_{*}$.
Applying theorem~\ref{theo:1d}  
once with $B=1$ and once with $B=\partial_{s}$,
and then applying proposition~\ref{prop:twofold}
yields a $C>0$ independent of $\Tt$ such that: 
\begin{equation} 
\label{eqK}
\|K\|\leq\sup(\|D\|,\|N\|)
\leq Ce^{\alpha L^{2}/T } = Ce^{\at_{*} \Lt^{2}/\Tt }
\ .
\end{equation}
We define $\vt\in C^{0}([0,\Tt], \mathcal{M}(]-\Lt,\Lt[))$
as the solution of 
\begin{equation*} 
\label{eqHeatvt}
\partial_{t}\vt - \partial_{s}^{2} \vt=0
\quad {\rm in}\ ]0,\Tt[\times ]-\Lt,\Lt[ ,\quad 
(\vt\res{s=-\Lt} ,\vt\res{s=\Lt})=b ,\quad
\vt\res{t=0} = \delta
\end{equation*}
where the control $b\in L^{2}(0,\Tt)^{2}$
is defined by $b(t)=0$ for $t\leq\eps \Tt$
and by $b(\eps\Tt+t')=K(\vt\res{t=\eps T})(t')$ for $t'\in ]0,T[$.
Note that $v_{0}=\vt\res{t=\eps T}$ is just 
the Dirac mass at the origin smoothed out 
by the homogeneous heat semigroup during a time $\eps \Tt$,
so that $v_{0}\in L^{2}(-L,L)$.
Moreover $\eps\Tt+T=\Tt$
and \linebreak $v(t,s)=\vt(\eps\Tt+t,s)$ is 
the solution of $(\ref{eqHeattwofold})$, 
so that $\vt\res{t=\Tt}=v\res{t=T}=0$.

To finish the proof that $\vt$
is a fundamental controlled solution
on \linebreak $]0,\Tt[ \times ]-\Lt,\Lt[$, 
we estimate its $L^{2}(]0,\Tt[\times ]-\Lt,\Lt[)$ norm 
which we abbreviate as $\|\vt\|_{\Tt,\Lt}$.
Setting $e_{j}(s)=\sin(j(s+\pi/2))\sqrt{2/\pi}$ defines 
an orthonormal basis $(e_{j})_{j\in \mathbb{N}^{*}}$
of $L^{2}(]-\Lt,\Lt[)$
such that $e_{j}$ is an eigenvector of $-\Delta_{s}$
with eigenvalue $j^{2}$.
In the weak topology, the Dirac mass can be decomposed in this basis as 
$\delta(s)=\sum_{j}e_{j}(0)e_{j}(s)$.
Note that 
%$e_{j}(0)=\sin(j\pi/2)\in \{0,1\}$, so that 
the sequence $(e_{j}(0))_{j\in \mathbb{N}^{*}}$ is bounded. % by $1$.
For $t\in]0,\Tt]$, we introduce the coordinates
$(\vt_{j}(t))_{j\in \mathbb{N}^{*}}$ of $\vt(t,\cdot) \in L^{2}(]-\Lt,\Lt[)$
in the Hilbert basis $(e_{j})_{j\in \mathbb{N}^{*}}$.
Using these coordinates 
and abbreviating the $L^{2}(]0,\Tt[)$ norm as 
$\|\cdot\|_{\Tt}$, 
the function $\vt$ and its norm write
\begin{equation}
\label{eqvt}
\vt(t,s)=\sum_{j}\vt_{j}(t)e_{j}(s)
\quad {\rm and } \quad
\|\vt\|_{\Tt,\Lt}^{2}
=\int_{0}^{\Tt}\sum_{j}|\vt_{j}(t)|^{2}\, dt
=\sum_{j} \|\vt_{j}\|_{\Tt}^{2} 
\ .
\end{equation}
As in~\cite{FR71}, these coordinates can be computed by 
$\vt_{j}(0)=e_{j}(0)$ and %the Duhamel formula
\begin{equation}
\label{eqvtj}
\vt_{j}(t) =  e^{-j^{2}t}\vt_{j}(0) + \int_{0}^{t}e^{-j^{2}(t-t')}
\left( e_{j}'(-\Lt)\vt(t',-\Lt) - e_{j}'(\Lt)\vt(t',\Lt) \right)
\, dt'
\ .
\end{equation}
Using Young's inequality to estimate the second term of the right hand side,
we have (since $\Tt< 4$, $|e_{j}'(\pm\Lt)|=|\vt_{j}(0)|=\sqrt{2/\pi}<1$)
\begin{eqnarray*}
\|\vt_{j}\|_{\Tt}
& \leq & |\vt_{j}(0)| \|e^{-j^{2}t} \|_{\Tt}
+ \|e^{-j^{2}t}\|_{L^{1}(]0,\Tt[)} 
\left(  |e_{j}'(-\Lt)| \|\vt(t',-\Lt)\|_{\Tt}
+ |e_{j}'(\Lt)| \|\vt(t',\Lt)\|_{\Tt} \right) \\
& \leq & \frac{4}{j}
\left( 1+\|\vt(t',-\Lt)\|_{\Tt}
+ \|\vt(t',\Lt)\|_{\Tt}\right)
\ .
\end{eqnarray*}
Hence equation $(\ref{eqvt})$ implies 
\begin{eqnarray*}
\|\vt\|_{\Tt,\Lt}^{2}
&\leq & 
\left( 1 +  \|\vt(t',-\Lt)\|_{\Tt}^{2}
+ \|\vt(t',\Lt)\|_{\Tt}^{2}\right)\sum_{j}\frac{4^{3}}{j^{2}}
%\\& = &
= 
\frac{4^{3}\pi^{2}}{6} 
\left(1 +\|Kv_{0}\|_{L^{2}(]0,\Tt[)}^{2}\right)
\ .
\end{eqnarray*}
% \begin{equation*}
% \|\vt\|_{\Tt,\Lt}^{2}
% \leq 
% \left(\sum_{j}\frac{1}{j^{2}}\right)
% \left( 1 + 4 \left( \|\vt(t',-\Lt)\|_{\Tt}^{2}
% + \|\vt(t',\Lt)\|_{\Tt}^{2}\right)\right)
% =
% \frac{\pi^{2}}{6} 
% \left(1 +4 \|Kv_{0}\|_{L^{2}(]0,\Tt[)}^{2}\right)
% \ .
% \end{equation*}
But there is an $A'>0$ independent of $\eps\Tt<1$ such that:
$$ 
 \|v_{0}\|_{L^{2}(]-\Lt,\Lt[)}^{2}
= \sum_{j}|\vt_{j}(\eps \Tt)|^{2}
\leq \sum_{j} e^{-2j^{2}\eps \Tt}
\leq \frac{A'}{\sqrt{\eps \Tt}} %A'/\sqrt{\eps \Tt}
\ .
$$
Hence equation $(\ref{eqK})$ yields a 
$C'>0$ independent of $\Tt$ such that: 
\begin{equation*}
\|\vt\|_{\Tt,\Lt}
\leq 
\frac{8\pi}{\sqrt{6}} 
\left(1 +2 \sqrt{\pi}\|K\|\|v_{0}\|_{\Lt}\right)
\leq \frac{C'}{\sqrt{\Tt}}e^{\at_{*} \Lt^{2}/\Tt }
\ .
\end{equation*}
Since $\at> \at_{*}$, 
there is an $\At>0$ independent of $\Tt$ such that:
\linebreak
$\displaystyle
\|\vt\|_{\Tt,\Lt}^{2} \leq \At e^{\at \Lt^{2}/\Tt }
$.
This completes the proof that $\vt$ is 
a fundamental controlled solution
on $]0,\Tt[\times ]-\Lt,\Lt[$
at cost $(\tilde{A},\tilde{\alpha})$.
\end{proof}

\subsection{The transmutation of waves into heat}

In a third step we perform a transmutation 
of an exact control for the wave equation 
into a null-control for the heat equation.
Our transmutation formula can be regarded  
as the analogue of Kannai's formula (\ref{eqKannai})
where the kernel $e^{-s^{2}/(4t)}/\sqrt{4\pi t}$,
which is the fundamental solution 
of the heat equation on the line,
is replaced by the fundamental controlled solution 
that we have constructed in the previous step. 
To ensure existence of an exact control for the wave equation
we use the geodesics condition of Bardos-Lebeau-Rauch
(already mentioned above theorem~\ref{theo:ub}):
\begin{thm}[\cite{BLR92}]
\label{theo:BLR}
If $L>L_{\Omega}$ then
for all  
$(w_{0}, w_{1})\in H^{1}_{0}(M)\times L^{2}(M)$
and all $(w_{2}, w_{3})\in H^{1}_{0}(M)\times L^{2}(M)$
there is a control function 
$f\in L^{2}(\mathbb{R}_{+}\times M)$
such that the solution 
$w\in C^{0}(\mathbb{R}_{+},H^{1}_{0}(M))\cap C^{1}(\mathbb{R}_{+},L^{2}(M))$
of the mixed Dirichlet-Cauchy problem 
(n.b.~the time variable is denoted by $s$ here):
\begin{equation} 
\label{eqWave}
\partial_{s}^{2}w - \Delta w={\bf 1}_{]0,L[\times \Omega} f 
\quad {\rm in}\ \mathbb{R}_{+}\times M, \quad 
w=0 \quad {\rm on}\ \mathbb{R}_{+}\times\partial M,
\end{equation}
with Cauchy data
$(w,\partial_{s}w)=(w_{0},w_{1})$ at $s=0$,
satisfies 
$(w,\partial_{s}w)=(w_{2},w_{3})$ at $s=L$.
Moreover, the operator 
$S_{W}: \left( H^{1}_{0}(M)\times L^{2}(M) \right)^{2}
\to L^{2}(\mathbb{R}_{+}\times M)$
defined by $S_{W}\left((w_{0},w_{1}),(w_{2},w_{3})\right)=f$ 
is continuous.
\end{thm}

We assume that theorem~\ref{theo:1d} holds for some rate $\alpha_{*}$.
Let $\alpha > \alpha_{*}$, 
\linebreak
$T\in]0,\inf(1,L_{\Omega}^{2})[$ 
and $L>L_{\Omega}$ be fixed from now on.
Let $A>0$ and 
\linebreak
$v\in L^{2}(]0,T[\times ]-L,L[)$ 
be the corresponding constant and fundamental controlled solution
given by proposition~\ref{prop:fundcontrol}.
We define $\vv\in L^{2}(\mathbb{R}^{2})$ 
as the extension of $v$ by zero,
i.e. $\vv(t,s)=v(t,s)$ on $]0,T[\times]-L,L[$
and $\vv$ is zero everywhere else.
It inherits from $v$ the following properties 
\begin{eqnarray}
\partial_{t}\vv - \partial_{s}^{2} \vv=0
\quad {\rm in }\  \mathcal{D}'(]0,+\infty[\times ]-L,L[)\ , 
\label{eqvv1} \\
%\vv\in C^{0}([0,+\infty[, \mathcal{M}(\mathbb{R})) %L^{2}(\mathbb{R}))
\vv\res{t=0}  =  \delta 
\quad {\rm and }\quad \vv\res{t=T}  =  0 \ ,
\label{eqvv2} \\
\|\vv\|_{L^{2}(]0,+\infty[\times\mathbb{R})}
\leq Ae^{\alpha L^{2}/T } \ .
\label{eqvv3}
\end{eqnarray}

Let $u_{0}\in H^{1}_{0}(M)$ be an initial data for the heat equation $(\ref{eqHeat})$.
Let $w$ and $f$
%$w\in C^{0}(\mathbb{R},H^{1}_{0}(M))$ and $f\in L^{2}(\mathbb{R}\times M)$ 
be the corresponding solution and control function 
for the wave equation obtained by applying theorem~\ref{theo:BLR}
with $w_{0}=u_{0}$ and $w_{1}=w_{2}=w_{3}=0$.
We define $\ww\in L^{2}(\mathbb{R};H^{1}_{0}(M))$ 
and $\ff\in L^{2}(\mathbb{R}\times M)$ 
as the extensions of $w$ and $f$ 
by reflection with respect to $s=0$,
i.e. $\ww(s,x)=w(s,x)=\ww(-s,x)$ 
and $\ff(s,x)=f(s,x)=\ff(-s,x)$ 
on $\mathbb{R}_{+}\times M$.
Since $w_{1}=0$, equation $(\ref{eqWave})$ imply 
\begin{equation} 
\label{eqWaveww}
\partial_{s}^{2}\ww - \Delta \ww={\bf 1}_{]-L,L[\times \Omega} \ff 
\quad {\rm in}\  \mathcal{D}'(\mathbb{R}\times M), \quad 
\ww=0 \quad {\rm on}\ \mathbb{R}\times\partial M,
\end{equation}

The main idea of our proof is to use $\vv$
as a kernel to transmute $\ww$ and $\ff$
into a solution $u$ and a control $g$ 
for $(\ref{eqHeat})$. 
Since $\vv\in L^{2}(\mathbb{R}^{2})$, 
$\ww\in L^{2}(\mathbb{R};H^{1}_{0}(M))$ 
and $\ff\in L^{2}(\mathbb{R}\times M)$,
the transmutation formulas 
\begin{equation} 
\label{eqtrans}
u(t,x)=\int_{\mathbb{R}}\vv(t,s)\ww(s,x)\, ds  
 \quad {\rm and } \quad 
g(t,x)=\int_{\mathbb{R}} %\int_{-L}^{L} 
\vv(t,s)\ff(s,x)\, ds 
\ ,
\end{equation}
define functions 
$u\in L^{2}(\mathbb{R};H^{1}_{0}(\times M))$ 
and $g\in L^{2}(\mathbb{R}\times M)$.
Since $\ww(s,x)=\partial_{s}\ww(s,x)=0$ for $|s|=L$,
equations $(\ref{eqWaveww})$ and 
$(\ref{eqvv1})$ imply 
\begin{equation}
  \label{equ1}
\partial_{t}u - \Delta u=\cd{1}_{]0,T[\times \Omega} g
\quad {\rm in }\  \mathcal{D}'(]0,+\infty[\times M) %\ , 
 \quad {\rm and } \quad
u=0 \quad {\rm on}\ ]0,T[\times\partial M,
\end{equation}
The property $(\ref{eqvv2})$ of $\vv$ implies 
\begin{equation}
  \label{equ2}
u\res{t=0}  =  u_{0} \quad {\rm and } \quad u\res{t=T}  =  0\ .
\end{equation}
%Since $S_{W}$ is continuous,
Setting $C=\sqrt{2}A\|S_{W}\|$,
Cauchy-Schwarz inequality with respect to $s$,
the estimate $(\ref{eqvv3})$ 
and
$\|\ff\|_{L^{2}(\mathbb{R}\times M)}^{2}
=2\|S_{W}\left((u_{0},0),(0,0))\right)\|_{L^{2}(\mathbb{R}_{+}\times M)}^{2}$
%$S_{W}\left((u_{0},0),(0,0))\right)=f$ 
imply
\begin{equation}
  \label{equ3}
\|g\|_{L^{2}(\mathbb{R}\times M)}
\leq \|\vv\|_{L^{2}(\mathbb{R}^{2})} \|\ff\|_{L^{2}(\mathbb{R}\times M)}
\leq Ce^{\alpha L^{2}/T }  \|u_{0}\|_{H^{1}_{0}(M)}\ .
\end{equation}

We have proved that 
for all $\alpha > \alpha_{*}$ there is a $C>0$
such that for all 
\linebreak
$u_{0}\in H^{1}_{0}(M)$, $T\in]0,\min\{1,L_{\Omega}^{2}\}[$ 
and $L>L_{\Omega}$, 
there is a control $g$ which solves 
the null-controllability problem %$(\ref{eqHeat})$
$(\ref{equ1})$, $(\ref{equ2})$, 
at a cost so estimated in $(\ref{equ3})$.
The same property holds for 
the space of data $L^{2}(M)$ instead of $H^{1}_{0}(M)$,
since $\|e^{\eps T\Delta}u_{0}\|_{H^{1}_{0}(M)}
\leq\|u_{0}\|_{L^{2}(M)}C_{0}/\sqrt{\eps T}$
with $\eps\in ]0,1[$ and 
$C_{0}=\|(1+\lambda)e^{-2\sqrt{\lambda}}\|_{L^{\infty}(\R)}^{1/2}$.
Therefore 
$\displaystyle
\limsup_{T\to 0} T \ln C_{T,\Omega} \leq \alpha
L^{2}$.
Letting $\alpha$ and $L$ tend respectively to $\alpha_{*}$ and $L_{\Omega}$
in this estimate completes the proof of $(\ref{eq:ub})$.

%\bibliographystyle{amsalpha}
%\bibliography{cost}

\providecommand{\bysame}{\leavevmode\hbox to3em{\hrulefill}\thinspace}
\providecommand{\MR}{\relax\ifhmode\unskip\space\fi MR }
\providecommand{\MRhref}[2]{%
  \href{http://www.ams.org/mathscinet-getitem?mr=#1}{#2}
}
\providecommand{\href}[2]{#2}

\end{document}